\documentclass{amsart}

\usepackage[utf8]{inputenc}
\usepackage{amsfonts}
\usepackage{mathtools}
\usepackage{amsmath}
\usepackage{xcolor}
\usepackage{amsthm}
\usepackage{pdflscape}
\usepackage{pgfplots}
\usepackage{mathrsfs}
\usepackage{hyperref}
\usepackage{tikz-cd}
\usepackage[normalem]{ulem}
\usepackage{amssymb}

\pgfplotsset{compat=1.18}

\newtheorem{thm}{Theorem}[section]
\newtheorem{cor}[thm]{Corollary}
\newtheorem{prop}[thm]{Proposition}
\newtheorem{lem}[thm]{Lemma}

\newtheorem{claim}[thm]{Claim}

\newtheorem{fact}[thm]{Fact}
\newtheorem*{thm*}{Theorem}

\theoremstyle{definition}
\newtheorem{defn}[thm]{Definition}

\theoremstyle{remark}
\newtheorem{rem}[thm]{Remark}

\newcommand{\partials}[2]{\frac{\partial #1}{\partial #2}}

\theoremstyle{plain}
\newtheorem{thmx}{Theorem}

\def\tp{\operatorname{tp}}

\def\forkindep{\mathrel{\raise0.2ex\hbox{\ooalign{\hidewidth$\vert$\hidewidth\cr\raise-0.9ex\hbox{$\smile$}}}}}

\begin{document}
\title[Algebraic independence for Lotka-Volterra systems]{Algebraic independence of solutions to multiple Lotka-Volterra systems}
\author{Yutong Duan}
\address{Yutong Duan\\
University of Illinois at Chicago\\
Department of Mathematics, Statistics and Computer Science\\
Science and Engineering Offices\\
Chicago, IL, 60607\\
United States}
\email{yduan24@uic.edu}

\author{Christine Eagles}
\address{Christine Eagles\\
University of Waterloo\\
Department of Pure Mathematics\\
Mathematics \& Computer\\
Waterloo, ON N2L 3G1\\
Canada}
\email{ceagles@uwaterloo.ca}

\author{L\'eo Jimenez}
\address{L\'eo Jimenez\\
The Ohio State University\\
Department of Mathematics\\
Math Tower\\
Columbus, OH 43210-1174\\
United States}
\email{jimenez.301@osu.edu}


\thanks{The first author was partially supported by the NSF grant DMS-2348885. The second author was supported by the Natural Sciences and Engineering Research Council of Canada (NSERC) [reference number CGS D - 588616 - 2024]}

\begin{abstract}\sloppy
    Consider some non-zero complex numbers $a_i, b_i, c_i, d_i$ with $1 \leq i \leq n$ and the associated classical Lotka-Volterra systems
    \[
    \begin{cases}
        X' = a_i XY + b_i X \\
        Y' = c_i XY + d_i Y \text{ .}
    \end{cases}
    \]
    \noindent We show that as long as $b_i \neq d_i$ for all $i$ and $\{ b_i, d_i\} \neq \{ b_j, d_j\}$ for $i \neq j$, any tuples $(x_1,y_1) , \cdots , (x_m,y_m)$ of pairwise distinct, non-degenerate solutions of these systems are algebraically independent over $\mathbb{C}$, meaning $\mathrm{trdeg}((x_1,y_1) , \cdots , (x_m,y_m)/\mathbb{C}) = 2m$. Our proof relies on extending recent work of Duan and Nagloo by showing strong minimality of these systems, as long as $b_i \neq d_i$. We also generalize a theorem of Brestovski which allows us to control algebraic relations using invariant volume forms. Finally, we completely classify all invariant algebraic curves in the non-strongly minimal, $b_i = d_i$ case by using machinery from geometric stability theory. 
\end{abstract}

\maketitle

\tableofcontents

\section{Introduction}

One of the main concerns of the model theory of differential fields is to develop a classification of solution sets of differential equations according to the complexity of their geometry. Of particular interest are methods to identify \emph{strongly minimal} sets. The solution set of a system of a differential equations of the form:

\[
\begin{cases}
    X_1' = f_1(X_1, \cdots , X_n) \\
    \vdots \\
    X_n' = f_n(X_1, \cdots , X_n)
\end{cases}
\tag{$\dagger$} \label{eq: sys-intro}
\]
for some $f_i \in F(X_1,\cdots, X_n)$, with $F < \mathbb{C}$, is said to be \emph{strongly minimal} if for any differential field extension $F < K$ and any solution $y_1, \cdots , y_n$ of the system, we have:
\[\mathrm{trdeg}(y_1, \cdots , y_n/K) = n \text { or $0$ .}\]

In practice, it is often required to add finitely many inequations to the system to obtain strong minimality, because some solutions may have transcendence degree smaller than $n$ over the base field. This is exemplified by the Lotka-Volterra systems below. 

To avoid this issue, it is convenient to consider the slightly weaker notion of \emph{minimality}. We call a solution $(y_1,\cdots , y_n)$ of (\ref{eq: sys-intro}) \emph{generic} if
\[\mathrm{trdeg}(y_1,\cdots , y_n /F) = n \text{ .}\]
The system is \emph{minimal} if for all generic solution $(y_1, \cdots , y_n)$ and all differential field extension $F < K$, we have
\[\mathrm{trdeg}(y_1, \cdots , y_n/K) = n \text{ or } 0 \text{ .}\]

Strongly minimal and minimal sets play a key role in the model-theoretic study of algebraic ordinary differential equations for two reasons. First, their geometry is well understood thanks to Zilber's dichotomy, proved by Hrushovski and Sokolovic \cite{hrushovskiminimal} (see also Pillay-Ziegler \cite{pillay2003jet}). This often yields some strong differential transcendence results, as we will see in this article. Second, through the concepts of semi-minimal analysis and domination-decomposition, they can be seen as the building blocks of solution sets of ordinary algebraic differential equations. 

Finding equations (or systems of equations) of order strictly greater than one with (strongly) minimal solution sets has proven to be difficult. Indeed, the list of equations which are known to be minimal is rather short, see for example the introduction of \cite{devilbiss2023generic}. The latest addition to this list is the set of non-degenerate solutions to the classical Lotka-Volterra system
\[
\begin{cases}\tag{$LV_{a,b,c,d}$}\label{LV}
    X' = aXY + bX \\
    Y' = cXY+dY \\
\end{cases}
\]
\noindent for some $a,b,c,d \in \mathbb{C} \setminus \{ 0\}$, which Duan and Nagloo, in \cite{duan2025algebraic}, showed to be strongly minimal if $\frac{b}{d} \not\in \mathbb{Q}$ (by non-degenerate, they mean $x \neq 0$ and $y \neq 0$). They also proved that the set of non-degenerate solutions to the 2d Lotka-Volterra system:
\[
\begin{cases}
X' = aXY + bX \\
Y' = cXY + dY^2\\
\end{cases}
\]
\noindent is strongly minimal, for any $a,b,c,d \in \mathbb{C} \setminus \{ 0 \}$. 

\sloppy Once strong minimality is known, a strengthening of the Zilber dichotomy for autonomous differential equations, together with a result of Brestovski \cite{brestovski1989algebraic}, yields a strong differential transcendence result, namely that any non-constant, non-degenerate, pairwise different solutions $(x_1,y_1), \cdots , (x_n,y_n)$ of (\ref{LV}) living in some differential field extension of $\mathbb{C}$ satisfy $\mathrm{trdeg}(x_1,y_1, \cdots , x_n , y_n /\mathbb{C}) = 2n$, which is the maximal possible transcendence degree.

Strong minimality is also equivalent to the non-existence of invariant curves, potentially over differential field extensions. More precisely, consider any planar system
\[
\begin{cases}\tag{$S$}\label{some system}
    X'  = f(X,Y) \\
    Y' = g(X,Y)
\end{cases}
\]
\noindent where $f,g \in \mathbb{C}(X,Y)$. This defines a rational vector field $s$, i.e. a section of the tangent bundle of $\mathbb{A}^2$, by letting $s(x,y) = \left( (x,y),(f(x,y),g(x,y)) \right)$. Associated to this rational vector field is its Lie derivative $\mathcal{L}_s$ which is defined, for any $(K,\delta)$ differential field  containing $\mathbb{C}$ and $P \in K[X,Y]$, as follows:
\[\mathcal{L}_s(P) = f(X,Y) \partials{P}{X} + g(X,Y) \partials{P}{Y} + P^{\delta}\]
\noindent where $P^{\delta}$ is the polynomial obtained by applying $\delta$ to the coefficients of $P$. A polynomial $P$ is \emph{invariant} if it divides $\mathcal{L}_s(P)$ in $K[X,Y]$, and an \emph{invariant curve} is the zero set of the ideal generated by a non-zero invariant polynomial. 

The connection between strong minimality and invariant curves can now be stated as follows: \emph{the system \ref{some system} is strongly minimal if and only if for any differential field extension $K$ of $\mathbb{C}$, there are no invariant curves over $K$.} Thus, showing strong minimality amounts to a full classification of invariant algebraic curves, over arbitrary differential field extensions. Note that if the system \ref{some system} has finitely many invariant curves, we obtain a strongly minimal set by removing these invariant curves. This is what happens for the Lotka-Volterra systems: the set of non-degenerate points is obtained by removing the invariant curves $X= 0$ and $Y = 0$. 

The present article further studies the algebraic transcendence properties of solutions of Lotka-Volterra systems through the realization of two main goals:
\begin{itemize}
    \item prove strong minimality of the classical Lotka-Volterra systems in the remaining open case of $\frac{b}{d} \in \mathbb{Q}$ but $b \neq d$,
    \item classify the possible algebraic relations between solutions of different Lotka-Volterra systems.
\end{itemize}

For the first point, the key new input is a theorem of Jaoui \cite[Theorem D]{jaoui2023density} which reduces the problem to the classification of complex invariant curves, already established by Duan and Nagloo \cite{duan2025algebraic}.

For the second point, note that there are obvious transformations between solution sets of Lotka-Volterra equations: the set $LV_{a,b,c,d}$ is in bijection with $LV_{1,b,1,d}$ through $(x,y) \to (cx,ay)$, and in bijection with $LV_{c,d,a,b}$ through $(x,y) \to (y,x)$. Using strong minimality and model-theoretic techniques, we prove that these are the only possible relations between solutions of classical Lotka-Volterra systems. This is our main result, which is Theorem \ref{theo: main-theo}: 

\begin{thmx}
    Let $LV_{a_i,b_i,c_i,d_i}, 1 \leq i \leq n$ be Lotka-Volterra systems 
    \[
    \begin{cases}
        X' = a_i XY + b_iX  \\
        Y' = c_i XY + d_iY 
    \end{cases}
    \]
    with $a_i,b_i,c_i,d_i \in \mathbb{C}$ and $b_i \neq d_i$ for all $i$. Let $F$ be a differential field extension of $\mathbb{C}$ and $(x_1,y_1), \cdots , (x_m,y_m) \in F$ be non-constant solutions of any of these systems, with $x_i,y_i \neq 0$ for all $i$. 
    
    If $\mathrm{trdeg}(x_1,y_1, \cdots , x_m,y_m/\mathbb{C}) < 2m$, then either:
    \begin{itemize}
        \item $(b_i,d_i) = (b_j,d_j)$ for some $i ,j$ and $(c_i x_i, a_i y_i) = (c_j x_j, a_j y_j)$,
        \item $(b_i,d_i) = (d_j,b_j)$ for some $i , j$ and $(c_ix_i,a_iy_i) = (a_jy_j, c_jx_j)$.
    \end{itemize}
\end{thmx}

This achieves a full classification of the algebraic relations between solutions of potentially different Lotka-Volterra systems, outside of the case $b =d$, which is of a very different nature. As was pointed out in \cite[Theorem 4.7]{eagles2024internality} (see also \cite[subsubsection 3.1.1]{duan2025algebraic}), the set of non-degenerate solutions of this system is not strongly minimal, since for any solution $(x,y)$ to $LV_{a,b,c,b}$, we have $(cx-ay)' = b(cx-ay)$, and so in particular $\mathrm{trdeg}(cx-ay/\mathbb{\mathbb{Q}}) = 1$, which implies $\mathrm{trdeg}(x,y/\mathbb{\mathbb{Q}}(cx-ay)) = 1$. However, using \cite[Theorem 4.7]{eagles2024internality}, as well as substantial model-theoretic machinery, we do achieve a complete classification of invariant algebraic curves, given in Theorem \ref{theo: LV b=d}:

\begin{thmx}
    The only irreducible invariant curves of $LV_{a,b,c,b}$ are the complex curves $X=0$, $Y=0$, $cX-aY=0$, as well as $cX-aY -z$ for $z$ any non-zero solution of $Z' = bZ$.
\end{thmx}

For the model theorist, we can restate this theorem as follows: the only definable subsets of the set of solutions of $LV_{a,b,c,b}$ are given by boolean combinations of singletons and the above irreducible invariant curves.

We now say a few words about the methods employed in the present work. The main obstacle to showing strong minimality is the need to check for the existence of invariant algebraic curves over arbitrary differential field extensions. In their work \cite{duan2025algebraic}, Duan and Nagloo use the \emph{degree of non-minimality}, an important invariant recently introduced by Freitag and Moosa \cite{freitag2023bounding}, to bypass this issue. 

It is a general model theory fact that to check strong minimality, it is enough to check it over differential fields generated by a fixed, finite number of independent, generic solutions of the system. The smallest such number is the degree of non-minimality. For systems of equations defined over $\mathbb{C}$, Freitag, Jaoui and Moosa showed in \cite{freitag2023degree} that the degree of nonminimality is 1, and this is a key ingredient in Duan and Nagloo's proof, as it allows them to work over field extensions generated by one solution of (\ref{LV}).

In this article, we use related, but different, methods to generalize their result. The central concept is of a proper almost fibration. A proper fibration of a planar system of the form (\ref{some system}), seen as a rational vector field $(\mathbb{A}^2,s)$, is a rational map $f$ to some algebraic curve $C$ equipped with a rational vector field $r$, such that
\begin{center}
\begin{tikzcd}
    T \mathbb{A}^2  \arrow[r, "df"] & TC \\
    \mathbb{A}^2 \arrow[r,"f"] \arrow[u, "s"] & C \arrow[u,"r"]
\end{tikzcd}
\end{center}
commutes. A proper almost fibration is a similar map, but now from a finite cover of $\mathbb{A}^2$ (also equipped with a rational vector field) to a curve $C$. 

It is a result of Jaoui and Moosa \cite[Theorem 6.8]{jaoui2022abelian} that the generic type of an autonomous system with no proper almost fibration is either minimal or a finite cover of a logarithmic-differential equation on a simple abelian variety $A$. In the case of planar systems (and more generally systems on affine space), the second possibility is ruled out using the observation that an algebraic group associated to the equation, its Galois group, must be linear (see \cite[Lemma 3.1]{eagles2024internality}), but also a definable subgroup of the simple abelian variety $A$. 

To show that the Lotka-Volterra system has no proper almost fibration, we use \cite[Theorem B]{jaoui2022generic} of Jaoui, which states that if there is a singular point not contained in any complex invariant curve, then there is no proper almost fibration.\footnote{Note that the use of Jaoui and Moosa's result is not, strictly speaking, necessary, as the criterion we use could be deduced entirely from work of Jaoui \cite{jaoui2023density}, where the abelian variety case is ruled out using different methods. However, the equivalence between no proper almost fibrations and strong minimality seemed worth pointing out.}

The main improvement over degree of nonminimality methods is that we do not need to consider \emph{any} differential field extension, thereby facilitating the computations by working exclusively over the complex numbers. Finally, to classify the possible relations between solutions, we generalize a theorem of Brestovski \cite[Theorem 2]{brestovski1989algebraic}, which was already used by Duan and Nagloo \cite{duan2025algebraic}. This theorem puts strong restrictions on the possible relations between solutions of very specific equations. We generalize it to constrain relations between multiple equations.

\medskip

Our article is organized as follows. Section \ref{sec: prelim} contains preliminaries on the model theory of differential fields, $D$-varieties and differential forms. In subsection \ref{subsec: brestovski}, we prove a slight generalization of a result of Brestovski \cite{brestovski1989algebraic} which controls the possible relations between solutions of strongly minimal equations of a specific form. Section \ref{sec: LV} is dedicated to the proof of our main theorems. In subsections \ref{subsec: complex inv curve}, we classify the complex invariant curves, and use that information to prove strong minimality of the set of non-degenerate solutions of $LV_{a,b,c,d}$, as long as $b \neq d$, and give a shorter proof of strong minimality of non-degenerate solutions of $LV_{a,b,c,d}^{2d}$. Using our generalization of Brestovski's theorem, we classify all possible algebraic relations between solutions to strongly minimal Lotka-Volterra systems in subsection \ref{subsec: ortho}. Finally, we classify the invariant curves of systems of the form $LV_{a,b,c,b}$ in subsection \ref{subsec: b is d}. 

\medskip

\textbf{Acknowledgements.}

The authors are grateful to R\'{e}mi Jaoui for pointing out how to use his work to prove strong minimality, i.e. using Lemma \ref{lem: complex-inv-curves} to prove Theorem \ref{theo: LV-sm}. They are also grateful to Ronnie Nagloo for his guidance, and putting the team together. Finally, they are thankful to an anonymous referee for their many useful suggestions.

\section{Preliminaries}\label{sec: prelim}

\subsection{Some model theory of differential fields}

In this article, we will make use of model-theoretic methods to study the Lotka-Volterra systems of equations. For that reason, we will work in the standard model theory setup of a sufficiently saturated model $(\mathcal{U}, \delta)$ of $\mathrm{DCF}_0$. Let us give some details, and refer the reader to \cite{marker2005model} for a more thorough introduction to the model theory of differential fields.

First, the set $\mathcal{U}$ is equipped with the structure of a field of characteristic zero, and $\delta$ is a derivation on $\mathcal{U}$. This means that $\delta$ is additive and satisfies the Leibniz rule that $\delta(ab) = a \delta(b) + \delta(a)b$ for all $a,b \in \mathcal{U}$. In this article, we will mostly denote this derivation by $\delta(a) = a'$.

Moreover, we require $\mathcal{U}$ to contain all possible solutions to differential equations and inequations. More precisely, let $K$ be a differential subfield of $\mathcal{U}$, meaning a subfield such that for all $a \in K$, we have $\delta(a) \in K$. Consider the following differential rings:
\begin{itemize}
    \item $K \{ X_1, \cdots, X_n \} = K[X_1 \cdots , X_n,X_1', \cdots , X_n',  X_1^{(2)}, \cdots, X_n^{(2)}, \cdots ]$ the ring of differential polynomials over $K$ in the variables $X_1, \cdots, X_n$, equipped with the natural derivation,
    \item $K \langle X_1, \cdots, X_n \rangle$ its field of fractions, equipped with the unique derivation extending the one on $K \{ X_1, \cdots, X_n \}$. We will also denote, for $a_1, \cdots , a_n$ some elements from $\mathcal{U}$, by $K\langle a_1, \cdots , a_n \rangle$ the differential field generated by the $a_i$ over $K$.
\end{itemize}

That $\mathcal{U}$ is differentially closed means that for any differential subfield $K$ and any finite system of differential equations and inequations in $K\langle X_1, \cdots , X_n \rangle$, if the system has a solution in a differential field extension $K<F$, then it already has a solution in $\mathcal{U}$. 

That $\mathcal{U}$ is saturated means that in fact, any system of differential equations and inequations, containing a number of equations strictly smaller than the cardinality of $\mathcal{U}$, if it has a solution in a differential field extension $K<F$, it has a solution in $\mathcal{U}$. We refer the reader to \cite{marker2005model} for a proof of the existence of a differentially closed field in any cardinality. In this document, we will assume that $\mathcal{U}$ has the cardinality of the continuum. 

Of particular interest is the field of constants of $\mathcal{U}$, which is
\[\mathcal{C} := \left\{ a \in \mathcal{U} : a' = 0 \right\} \text{ .}\]
This field can be showed to be algebraically closed. Because we assume that $\mathcal{U}$ is the cardinality of the continuum, we can also assume that $\mathcal{C} = \mathbb{C}$, the field of complex numbers, and will do so in the rest of this article. 

We now recall the key model-theoretic notion of \emph{definable sets} in this context. By a definable set, we simply mean a boolean combination of solution sets, in $\mathcal{U}$, of differential equations with parameters in a subfield of $\mathcal{U}$. Note that the general notion of definable set is more complex, as it allows for quantification, or equivalently, states that a coordinate projection of a definable set is definable. That this coincides with the previous definition is a consequence of quantifier elimination for $\mathrm{DCF}_0$. 

The \emph{order} of a definable set $S$ over some field $K$ is
\[\mathrm{ord}(S) := \sup \{ \mathrm{trdeg}(K\langle a \rangle /K) : a \in S \} \]
\noindent and any element of $S$ realizing this maximum is called \emph{generic}, if the order is finite\footnote{One can also define order and genericity in the infinite case using the \emph{Kolchin polynomial} but we will not need this here.}. 

The notion of \emph{type} will also be crucial.

\begin{defn}
    Let $K < \mathcal{U}$ be a differential subfield. A \emph{type} $p(x_1, \cdots , x_n)$ over $K$ is a maximal set of differential equations and inequations in the variables $x_1, \cdots x_n$, with coefficients in $K$, having a solution in some differential field extension $K<L$.

    The set of types in $n$ variables over $K$ is denoted $S_n(K)$, and we write $S(K)$ when the number of variables is irrelevant or implicit.
\end{defn}

Types are often denoted $p,q,r,s, \cdots$. If $p \in S(K)$ is a type, a \emph{realization} of $p$ is a tuple of $\mathcal{U}$ satisfying all equations and inequations in $p$. 

Any tuple in $\mathcal{U}$ gives rise to a type: if $a_1, \cdots , a_n$ are elements of $\mathcal{U}$, the type of $a_1, \cdots , a_n$ over $K$ is the set of all differential equations and inequations, with coefficients in $K$, satisfied by $a_1,\cdots ,a_n$. It is denoted $\tp(a_1, \cdots , a_n/K)$.

If $K < \mathcal{U}$ is countable, then any type $p \in S(K)$ has a realization in $\mathcal{U}$: this is a consequence of saturation of $\mathcal{U}$. Because of this, we will often work over a countable subfield of $\mathcal{U}$. 
 
A word of warning for the reader unfamiliar with model theory: in this article, we will sometimes use the standard model theoretic convention of calling tuples elements when the distinction between the two is irrelevant. 

A useful invariant of a type is its \emph{order}:

\begin{defn}
      The \emph{order} of a type over $K$ is defined as the transcendence degree over $K$ of the differential field generated, over $K$, by any of its realizations. 
\end{defn}

The following definition is crucial in the model-theoretic study of $\mathrm{DCF}_0$ (as well as a much broader class of structures): 

\begin{defn}
    Given a finite order type $p = \tp(a/F)$ and a differential subfield $K<F$, we say that $p$ \emph{does not fork over $K$} if the order of $p$ is equal to the order of $\tp(a/K)$. In other words, the type $p$ does not fork over $K$ if and only if
    \[\mathrm{trdeg}(K\langle a\rangle/K)=\mathrm{trdeg}(F\langle a\rangle/F) \text{ .}\] 
    \noindent Note that this does not depend on which $a \models p$ we pick.

    We will also say that $a$ does not fork with $F$ over $K$, or that $a$ is independent from $F$ over $K$, and write $a \forkindep_K F$. If $b$ is another tuple, we write $a \forkindep_K b$ for $a \forkindep_K K\langle b \rangle$.

    We say that $\{a_1,\dots,a_n\}$ is an \emph{independent set over $K$} if for all $i=1,\dots,n$, $a_i\forkindep_K a_1,...,a_{i-1},a_{i+1},...,a_n$.
\end{defn}

This notion of independence has many useful properties, and agrees with Shelah's notion of forking independence in general stable theories, see Sections 7.2 and 7.3 of \cite{Tent_Ziegler_2012}. 

When $K$ is algebraically closed, for any $p \in S(K)$ and any field extension $K<F$, any two $a,b$ realizing $p$ and independent of $F$ over $K$ have the same type over $K$. We denote this type by $p{\upharpoonright}_F$. In general, we have the following definition:

\begin{defn}
    A type $p \in S(K)$ is called \emph{stationary} if for any field extension $K<F$, the type of a realization of $p$ independent from $F$ over $K$ is unique. In that case, we denote the type of a non-forking extension to $F$ by $p {\upharpoonright} _F$.
\end{defn}

This allows us to define the Morley product of types, which will only be used in subsection \ref{subsec: b is d}. Let $p, q\in S(K)$ be stationary types. Their \emph{Morley product} $p\otimes q$ is defined as $\tp(ab/K)$ where $a\models p$ and $b\models q{\upharpoonright}_{K\langle b\rangle}$. By $p^{(n)}$ we mean $\tp(a_1,\dots,a_n/K)$ such that for all $i=1,\dots,n$, $a_i\models p{\upharpoonright}_{K\langle a_1,...,a_{i-1}\rangle}$.

\medskip

We now come to the central notion of strong minimality. 

\begin{defn}
    Let $S$ be a set definable over some differential field $K$. It is \emph{strongly minimal} if it is infinite and for any other definable set $R$ (over any set of parameters), either $S \cap R$ or $S \setminus R$ is finite.
\end{defn}

Note that the definable set $S$ need not be a subset of $\mathcal{U}$: it could also be a subset of $\mathcal{U}^n$, for some $n>1$. We consider a concrete example below.

If $S$ is a definable set given by the solutions of a system:
\[
\begin{cases}
    X_1' = f_1(X_1, \cdots, X_n) \\
    \vdots \\
    X_n ' = f_n(X_1, \cdots, X_n)
\end{cases}
\]
\noindent with the $f_i \in K(X_1, \cdots , X_n)$, strong minimality of $S$ is equivalent to the condition that for any $a \in S$ and any differential field extension $K<F$, we have:
    \[\mathrm{trdeg}(F \langle a \rangle/F) = 0 \text{ or $n$.}\]

We now give a brief explanation of this fact. Suppose $S$ is strongly minimal, and let $a= (a_1, \cdots , a_n) \in S$. If $\mathrm{trdeg}(F\langle a \rangle/F) < n$ for some field extension $F$, then there must be some $P \in F[X_1, \cdots , X_n] \setminus \{ 0 \}$ such that $P(a_1, \cdots , a_n) = 0$. Consider the definable set $R = \{ (x_1, \cdots , x_n) \in \mathcal{U} : P(x_1, \cdots , x_n) = 0\}$, then either $S \cap R$ is finite or $S \setminus R$ is finite. But if $S \setminus R$ is finite, since there are infinitely many elements of $S$ of maximal transcendence degree, we conclude that $S \cap R$ must contain some $b_1, \cdots , b_n$ with $\mathrm{trdeg}(F(b_1, \cdots , b_n)/F) = n$. This implies $P = 0$, a contradiction. So $S \cap R$ is finite. By quantifier elimination and differential closedness of $\mathcal{U}$, this implies $\mathrm{trdeg}(F\langle a_1, \cdots , a_n \rangle/F) = 0$. The other direction involves similar considerations.

In subsection \ref{subsec: D-var}, we will give a third equivalence condition in terms of non-existence of invariant curves.  

We can also define strong minimality for types:

\begin{defn}
    We say that a type $p \in S(K)$ is \emph{strongly minimal} if some (any) of its realization is not in $\mathrm{K}^{\mathrm{alg}}$, and it contains a formula defining a strongly minimal set. 
\end{defn}

For types, the slightly weaker notion of \emph{minimality} if often easier to work with:

\begin{defn}
    We say that a type $p \in S(K)$ is \emph{minimal} if for any $a$ realizing $p$ and extension $K<F$, either $a \forkindep_K F$ or $a \in F^{\mathrm{alg}}$. 
\end{defn}

Note that strong minimality and minimality for types can be defined in any first-order theory. In general strong minimality is strictly stronger, but to our knowledge, no type in $\mathrm{DCF}_0$ that is minimal but not strongly minimal has been identified.

Equivalently, a type $p \in S(K)$ is minimal if and only if for any $a$ realizing $p$ and any differential field extension $K<F$, we have
\[\mathrm{trdeg}(F \langle a \rangle/F) = \mathrm{trdeg}(K \langle a \rangle/K) \text{ or $0$.}\]
\noindent As was pointed out in the introduction, strongly minimal sets and minimal types can be classified according to the complexity of their internal geometry. An example of particular interest to us are strongly minimal sets with a very rudimentary structure:

\begin{defn}
    Let $S$ be a strongly minimal set defined over $K$. We say that $S$ is \emph{disintegrated} (sometimes also called geometrically trivial) if for any distinct $x_1, \cdots , x_m \in S \setminus K^{\mathrm{alg}}$, if the $x_i$ and their derivatives are not algebraically independent over $K$, then there are $i \neq j$ such that already $x_i, x_j$ and their derivatives are not algebraically independent over $K$. In terms of transcendence degree, if $n$ is the order of $S$, then this means:
    \[\mathrm{trdeg}(K\langle x_1, \cdots , x_m \rangle/K) < mn \Rightarrow \mathrm{trdeg}(K \langle x_i, x_j \rangle/K) < 2n \text{ for some $i \neq j$.}\]

    The set $S$ is \emph{totally disintegrated} if for any distinct $x_1, \cdots , x_m \in S \setminus K^{\mathrm{alg}}$, the $x_i$ and their derivatives are algebraically independent over $K$. In terms of transcendence degree:
    \[\mathrm{trdeg}(K \langle x_1, \cdots, x_m \rangle/K) = nm \text{ for all distinct $x_1, \cdots , x_n \in S\setminus K^{\mathrm{alg}}$.}\]
\end{defn}

It is good to keep in mind the following intuition: in a disintegrated strongly minimal set, relations between its elements are controlled by relations between pairs of elements. If the set is totally disintegrated, there is no relation between its elements besides equality.

We can make similar definitions for (strongly) minimal types, replacing elements of a definable set by realizations of a type. In other words:

\begin{defn}
    Let $p \in S(K)$ be a minimal type. It is \emph{disintegrated} if for any distinct realizations $x_1, \cdots , x_n$ of $p$, we have
    \[\mathrm{trdeg}(K\langle x_1, \cdots , x_m \rangle/K) < mn \Rightarrow \mathrm{trdeg}(K \langle x_i, x_j \rangle/K) < 2n \text{ for some $i \neq j$.}\]
    \noindent It is \emph{totally disintegrated} if for any distinct realizations $x_1, \cdots , x_n$, we have
    \[\mathrm{trdeg}(K \langle x_1, \cdots, x_m \rangle/K) = nm \text{ .}\]
\end{defn}

Here we see why working with types instead of definable sets is convenient: we do not have to remove algebraic solutions.

As a consequence of the classification of the non-disintegrated strongly minimal sets in $\mathrm{DCF}_0$, due to Hrushovski and Sokolovic \cite{hrushovskiminimal}, we can show that most strongly minimal sets over constant parameters are disintegrated (see for example \cite[Proposition 5.8]{casale2020ax}):

\begin{fact}
Let $S$ be a strongly minimal set of order $>1$ and suppose that $S$ is defined over $\mathbb{C}$. Then $S$ is disintegrated.
\end{fact}

Again, we can write the same statements for types:

\begin{fact}\label{fact: order-2-sm-disin}
    Let $p \in S(K)$ be a type of order strictly greater than one, for some $K < \mathbb{C}$. If $p$ is (strongly) minimal, then is it disintegrated.
\end{fact}

Using this fact, we can deduce, from strong minimality, some very strong differential transcendence results, since possible relations between realizations of $p$ will be controlled by relations between pairs of solutions. 

Total disintegration says that there is as little as possible relations between solutions of a strongly minimal equation. The equivalent concept for pairs of equations is (weak) \emph{orthogonality}, which is better defined for types:

\begin{defn}
    \begin{itemize}
        \item[]
        \item Let $p,q\in S(K)$. We say \emph{$p$ is weakly orthogonal to $q$} if whenever $a$ realizes $p$ and $b$ realizes $q$, then $a\forkindep_K b$. We denote this by $p\perp^w q$. 
        \item Let $p\in S(F)$ and $q\in S(K)$ be stationary types. Then \emph{$p$ is orthogonal to $q$}, denoted $p\perp q$, if for every differential field $L$ containing $K$ and $F$, we have $p{\upharpoonright}_L\perp^wq{\upharpoonright}_L$. 
    \end{itemize}
\end{defn}

So intuitively, two types $p$ and $q$ are weakly orthogonal if there is no relations between any two pairs of their realizations. They are orthogonal if this stays true after taking a non-forking extension. Orthogonality clearly implies weak orthogonality, but the converse is false (see the discussion below Definition \ref{def: weak-ortho-to-C} for a related example).

Again, disintegratedness is useful here, by the following fact:
\begin{fact}\cite[Corollary 2.5.5]{pillay1996geometric} \label{weak non-orthogonality}
    If $p\in S(F)$ and $q\in S(K)$ are two minimal and disintegrated types which are weakly orthogonal, then they are orthogonal.
\end{fact}

Often of particular interest is whether a type $p$ has no definable relations with the field of constants. 

\begin{defn}\label{def: weak-ortho-to-C}
    Let $p\in S(K)$. Then $p$ is \emph{weakly orthogonal} to $\mathbb{C}$ (or to the constants) if for all $a$ realizing $p$ and every finite tuple $c \subseteq \mathbb{C}$, we have $a\forkindep_K c$. 
    
    If $p$ is stationary, we say $p$ is \emph{orthogonal to the constants} if for every differential field extension $K<L$, the type $p{\upharpoonright}_L$ is weakly orthogonal to $\mathbb{C}$.
\end{defn}

Again, orthogonality is strictly stronger than weak orthogonality. Consider for example the generic type $p$ of the definable set $y'= y$ over $\mathbb{Q}$. One can show that for any realization $a$ and tuple $c$ of complex numbers, we have $a \forkindep_{\mathbb{Q}} c$. However, if $b$ is another realization of $p$ with $a \forkindep_{\mathbb{Q}} b$, then $\left( \frac{a}{b}\right)' \in \mathbb{C}$, and thus $p {\upharpoonright} _{\mathbb{Q}(b)}$ is not weakly orthogonal to $\mathbb{C}$.

Important to us will be the fact that any minimal type of order strictly larger than one is orthogonal to $\mathbb{C}$:

\begin{fact}\label{fact: mini-implies-ortho}
    Any minimal type $p \in S(K)$ of order strictly larger than one is orthogonal to the constants.
\end{fact}

An opposite notion to orthogonality to $\mathbb{C}$ is \emph{(almost) internality} to $\mathbb{C}$. 

\begin{defn}
    A type $p \in S(K)$ is internal to $\mathbb{C}$ if it is stationary and there is a differential field extension $K<L$ and some $a$ realizing $p$ such that both:
    \begin{itemize}
        \item $a \forkindep_K L$,
        \item $a \in L(c_1, \cdots , c_k)$ for some $c_1 , \cdots , c_k \in \mathbb{C}$.
    \end{itemize}
    If we have $a \in L(c_1, \cdots , c_k)^{\mathrm{alg}}$ instead, we say $p$ is \emph{almost} $\mathbb{C}$-internal.
\end{defn}

For example, it is not difficult to see that if $a$ is a solution of a \emph{linear} differential equation over some algebraically closed differential field $K$, then $\tp(a/K)$ is $\mathbb{C}$-internal. Indeed, the set of solutions to $L(x) = 0$ is a finite-dimensional $\mathbb{C}$-vector space, and one can pick $L$ to be the extension generated by a basis $B$ satisfying $a \forkindep_K K\langle B \rangle$. That such an independent basis exists is a consequence of the behavior of forking in $\mathrm{DCF}_0$, and more generally in stable theories. 

We will briefly make use of the \emph{binding group} of a $\mathbb{C}$-internal type $p \in S(K)$ in the proof of Theorem \ref{theo: min-crit} below. We briefly recall a few important facts. The binding group of $p$ is a definable group acting faithfully on the set of realizations of $p$. By definable group, we mean that its set of elements, the graph of its multiplication and inverse functions are definable sets. Moreover, its action is isomorphic to the action of $\mathrm{Aut}_K(p/\mathbb{C})$, the group of differential field automorphisms of the ambient field, fixing $\mathbb{C}$ and $K$ pointwise, on $p(\mathcal{U})$. Finally, this binding group is definably isomorphic to the $\mathbb{C}$-points of an algebraic group. Note that in the case of the type of a solution to a linear differential equation, this algebraic group is in fact the differential Galois group familiar to differential algebraists.

A type  might not be almost internal to the constants but instead `built up' from types almost internal to the constants. More precisely, we say $p\in S(K)$ is $\mathbb{C}$-\emph{analyzable} if there are a realization $a$ of $p$ and a sequence $a = a_n, a_{n-1}, \cdots, a_1 $ such that:
\begin{itemize}
    \item the type $\tp(a_i/K\langle a_{i-1} \rangle)$ is almost $\mathbb{C}$-internal for all $i > 1$, and $\tp(a_1/K)$ is almost $\mathbb{C}$-internal,
    \item $a_{i-1} \in K\langle a_i \rangle$ for all $i>1$.
\end{itemize}

We say a type is \emph{$n$-step analyzable} if there is an analysis of length $n$. Analyzability will only play a role in Subsection \ref{subsec: b is d}.

\subsection{\texorpdfstring{$D$}{D}-varieties and rational forms}\label{subsec: D-var}

We now recall the correspondence between polynomial and rational (twisted) vector fields, $D$-varieties, and derivations, and refer the reader to \cite{moosa2022six} for more details. 

Let $K$ be any differential field, and $V \subset \mathcal{U}^n$ some irreducible affine variety defined over $K$. The \emph{prolongation} of $V$ is the variety $\tau V \subset K^{2n}$ defined by the following equations:
\[
\begin{cases}
    P(X_1, \cdots, X_n) = 0 \\
    P^{\delta}(X_1, \cdots , X_n) + \sum\limits_{i=1}^n \partials{P}{X_i} U_i = 0
\end{cases}
\]
\noindent for each $P \in I(V)$, where $P^{\delta}$ is the result of applying the ambient derivation to the coefficients of $P$. Note that if $K < \mathbb{C}$, we recover the classical tangent bundle $TV$. We get a surjective morphism $\pi : \tau V \rightarrow V$.

Taking $V$ to $\tau V$ defines a functor as we can apply it to the graphs of morphisms. Moreover, we can, by patching, define $\tau V$ for any algebraic variety.

We now define an algebraic object central to our methods:

\begin{defn}
    A \emph{$D$-variety} is an algebraic variety $V$ equipped with a regular section $(V,s)$ of its prolongation. 
    
    If $s$ is a rational section of its prolongation instead, we call $(V,s)$ a rational $D$-variety.
\end{defn}

Note that if $V$ is defined over some $K < \mathbb{C}$, this is simply a polynomial (resp. rational) vector field, i.e. a regular (resp. rational) section of the tangent bundle.

It is straightforward to check the following crucial result:

\begin{fact}
    The set of regular (resp. rational) sections of the prolongation $\tau V$ of a variety $V$ is in one-to-one correspondence with on $K[V]$ (resp. $K(V)$) extending the derivation on $K$.
\end{fact}

We refer the reader to \cite{moosa2022six} for a rigorous explanation of that fact. 

To any $D$-variety is canonically associated a definable set, its set of $D$-points, or sharp points:

\begin{defn}
    Let $(V,s)$ over $K$ be a $D$-variety with $V\subseteq\mathbb{A}^n$ and $s(\bar{x})=(\bar{x},v_1(\bar{x}),...,v_n(\bar{x}))$. 
    The set of \emph{$D$}-points, or \emph{sharp points}, of $(V,s)$ is the set of $(a_1,...,a_n)\in V(\mathcal{U})$ such that $\delta (a_i)=v_i(a_1, \cdots , a_n)$. It is denoted by $(V,s)^{\sharp}$. 
\end{defn}

We can also define the set of sharp points of a rational $D$-variety $(V,s)$, by adding the condition that the section $s$ is well-defined.

As an important special case, consider a variety $V$ defined over some $K < \mathbb{C}$. Then $\tau V = TV$, and $(a,0) \in TV$ for any $a \in V$. Thus, the map $s$ sending any $a \in V$ to $(a,0)$ is a section of $TV$, denoted $0$, and $(V,0)$ is a $D$-variety. The sharp points $(V,0)^{\sharp}$ corresponds to $V(\mathbb{C})$, the constant points of $V$.

As another example, consider a system of the form
\[
\begin{cases}
    y_1' = f_1(y_1, \cdots , y_n) \\
    \vdots \\
    y_n' = f_n(y_1, \cdots , y_n)
\end{cases}
\]
\noindent with $f_i \in K(X_1, \cdots , X_n)$. It corresponds to the rational $D$-variety $(\mathbb{A}^n, s)$, with 
\begin{align*}
    s: \mathbb{A}^n & \rightarrow \tau \mathbb{A}^n \\
    (y_1, \cdots , y_n) & \rightarrow ((y_1,\cdots, y_n), (f_1(y_1,\cdots , f_n), \cdots , f_n(y_1, \cdots , y_n))
\end{align*}
\noindent and the sharp points of $(\mathbb{A}^n,s)$ are exactly the solutions of the system.

The formalism of $D$-varieties makes it easy to define generic points:

\begin{defn}
    The \emph{generic type} over a differential field $K$ of a (rational) $D$-variety defined over $K$ is the type of some $D$-point $a\in (V,s)^{\sharp}$ such that $V$ is the Zariski locus of $a$ over $K$. A \emph{generic point} of $(V,s)$ is a realization of its generic type.
\end{defn}

Finite order types are entirely captured by rational $D$-varieties: 

\begin{fact}\label{fact: d-var-capture-types}
    Let $p \in S(K)$ be a type of finite order. There is a rational $D$-variety $(V,s)$, defined over $K$, such that $p$ is interdefinable with the generic type of $(V,s)$.
\end{fact}

We refer the reader to \cite{moosa2022six} for a proof of this fact. 

We now give an example of this phenomenon that will be relevant later. Consider $K<\mathbb{C}$, some $f,g \in K(X_0,X_1)$, and the equation:
\[(E): x''f(x,x')+g(x,x') = 0 \text{ .}\]
This definable set of solutions of $(E)$ corresponds exactly, under the definable map $x \rightarrow (x,x')$, to the $D$-points of the variety given by the section
\[s: (x_0,x_1) \rightarrow \left((x_0,x_1), \left(x_1,\frac{-g(x_0,x_1)}{f(x_0,x_1)}\right)\right)\]
of the tangent bundle of $\mathbb{A}^2$, the affine plane. If we now consider the generic type $p$ of the set of solutions of $(E)$, it corresponds to solutions $x$ such that $(x,x')$ satisfy no algebraic relation. The map $x \rightarrow (x,x')$ induces a bijection between these and generic points of $(\mathbb{A}^2,s)$.

If we wanted to obtain a definable bijection between $p$ and the generic points of a $D$-variety, we could do it by introducing the extra variable $y$, and considering the section of the tangent bundle of $V = \{ (x_0,x_1,y) : yf(x_0,x_1) = 1\}$ given by $(x_1, -g(x_0,x_1), y')$, noting that $y'$ can be expressed as a polynomial in $x_0,x_1$ and $y$.

We can define morphisms of $D$-varieties:

\begin{defn}
    Given $(V,s)$ and $(W,r)$ some (rational) $D$-varieties over $K$, a morphism from $(V,s)$ to $(W,r)$ is a (rational) regular map $f$ such that the following
    \[
    \begin{tikzcd}
        \tau V \arrow[r, "\tau f"] & \tau W \\
        V \arrow[r, "f"] \arrow[u, "s"] & W \arrow[u, "r"]
    \end{tikzcd}
    \]
    commutes. 
    
    In particular, a $D$-subvariety (also called an invariant subvariety) of $(V,s)$ is a subvariety $W \subset V$ over $K$ such that the inclusion map is a morphism of $D$-varieties. Equivalently, the subvariety $W$ is a $D$-subvariety of $s$ if and only if $s(W) \subset \tau W$.
\end{defn}

We encourage the reader unfamiliar with the formalism of $D$-varieties to try and prove the following fact: sharp points of a $D$-variety $(V,s)$ are exactly points $a \in V$ such that $x=a$ is a $D$-subvariety.

There is a natural connection between strong minimality and invariant subvarieties. As remarked previously, the set of sharp points of a $D$-variety form a definable set, and we may ask if it is strongly minimal. We have the following, which is proved in \cite[Section 2.2]{nagloo2017algebraic}:

\begin{fact}
    The set of sharp points of a $D$-variety $(V,s)$ is strongly minimal if and only if it has no proper positive dimensional $D$-subvarieties, over any field of definition.
\end{fact}

Reasoning along these lines, we can obtain a stronger result: if a $D$-variety $(V,s)$ defined over $K$ has only finitely many proper positive dimensional $D$-subvarieties $W_1, \cdots , W_m$, the definable set obtained by removing their sharp points, i.e. $(V,s)^{\sharp} \setminus \left( \bigcup\limits_{i=1}^m (W_i,s)^{\mathrm{\sharp}}\right)$, is strongly minimal. This will be the case for the classical Lotka-Volterra systems, where we remove $X = 0$ and $Y = 0$. Bottom line, even if the definable set of sharp points of a $D$-variety is not strongly minimal, it may become strongly minimal after removing finitely many $D$-subvarieties. 

Establishing strong minimality of the generic type of a $D$-variety $(V,s)$ therefore can be obtained by classifying its positive dimensional proper invariant subvarieties. This problem is generally of interest, and we also solve it for a non-minimal Lotka-Volterra system in Subsection \ref{subsec: b is d}.

Because we will mostly work with Lotka-Volterra systems, which are $D$-varieties on $\mathbb{A}^2$, dimension one $D$-subvarieties will be of particular interest. We will call them \emph{invariant curves}. To prove that the generic type of a $D$-variety on $\mathbb{A}^2$ is strongly minimal, it is enough to fully classify its invariant curves.

We can characterize $D$-subvarieties using the derivative associated to a $D$-variety. As a special case, consider a $D$-variety $(\mathbb{A}^2,s)$ given by a system of the form
\[
\begin{cases}
    X' = f(X,Y) \\
    Y' = g(X,Y)
\end{cases}
\]
\noindent with $f,g \in \mathbb{C}[X,Y]$, as well as some $P \in \mathbb{C}[X,Y]$. Let $D_s$ be the associated derivation. Then the algebraic curve described by $P(X,Y) = 0$ is an invariant curve of $(\mathbb{A}^2, s)$ if and only if $P$ divides $D_s(P)$.

\subsection{Some differential geometry tools}

The previously established correspondence between differential equations and rational vector fields will allow us to use tools from differential geometry, in particular Lie derivatives and differential forms. For the comfort of the reader, we will give proofs of some of the results we state, but note that nothing in this section is new. For more on differential forms, we refer the reader to any of the many books on the subject, for example \cite{lang2012fundamentals}.

The first step to define rational differential forms will be to dualize the space of rational sections.

Fix an $n$-dimensional irreducible affine variety $V$ defined over $\mathbb{C}$ and let $\mathbb{C}(V)$ be its function field. As discussed previously, we identify the $\mathbb{C}(V)$-vector spaces of derivations $\mathrm{Der}(\mathbb{C}(V)/\mathbb{C})$ and of rational vector fields on $V$. 

Locally, given $p \in V$, a section $s$ gives us $s(p)$, an element of $T_pV$, the tangent space of $V$ at $p$. Since $T_pV$ is a $\mathbb{C}$-vector space, it has a dual, which is the cotangent space $T_p^*V$ of $V$ at $p$.

Going back to the global picture, sections of the tangent bundle form the $F$-vector space of rational vector fields on $V$, which has a dual $\mathbb{C}(V)$-vector space, denoted $\Omega^1_V(\mathbb{C}(V)) = \Omega^1_V$. If $\omega \in \Omega^1_V$, then for any $p \in V$, its evaluation $\omega(p)$ is an element of $T_p^*V$. 

\begin{defn}
    The dual vector space of the $\mathbb{C}(V)$-vector space $\mathrm{Der}(\mathbb{C}(V)/\mathbb{C})$ is called the space of \emph{rational $1$-forms}, and denoted $\Omega_V^1$.
\end{defn}

Given a fixed $f \in \mathbb{C}(V)$, evaluating derivations $D \in \mathrm{Der}(\mathbb{C}(V)/\mathbb{C})$ at $f$ is a $\mathbb{C}(V)$-linear map from $\mathrm{Der}(\mathbb{C}(V)/\mathbb{C})$ to $\mathbb{C}(V)$. In other words, we naturally obtain the \emph{exterior derivative}:

\begin{defn}
    The exterior derivative $d : \mathbb{C}(V) \rightarrow \Omega^1_V$ is defined as the map
    \begin{align*}
        d : \mathbb{C}(V) & \rightarrow \Omega^1_V \\
        f & \rightarrow (D \rightarrow D(f))
    \end{align*}
    It is a $\mathbb{C}$-linear map, and $d(fg) = df \times g + f \times dg$ for all $f,g \in F$.
\end{defn}

We point out that we can identify bases of the vector spaces of derivations and $1$-forms. Indeed, let $u_1,\cdots ,u_n$ be a transcendence basis of $\mathbb{C}(V)$ over $\mathbb{C}$. We define derivations $\partials{}{u_i}$ by setting $\partials{}{u_i}(u_j) = 1$ if $j = i$ and $0$ otherwise. It is a standard differential algebra fact that this extends to a unique derivation on $\mathbb{C}(V)$. By definition, we have that $du_i(\partials{}{u_j}) = \partials{}{u_j}(u_i)$, and therefore $(du_i)_{i = 1 \cdots n}$ is the dual basis of $\Omega_V^1$. In particular $\Omega_V^1$ is an $\mathbb{C}(V)$-vector space of dimension $\dim(V)$. 

One useful feature of the exterior derivative is that is allows us to prove certain functions are constant:

\begin{prop}\label{prop: df-zero-implies-cst}
    Let $f \in \mathbb{C}(V)$. Then $f \in \mathbb{C}$ if and only if $df = 0$.
\end{prop}

\begin{proof}
    The left to right direction is immediate, since $\Omega_V^1$ is the $\mathbb{C}$-vector spaces of $\mathbb{C}$-linear derivations on $\mathbb{C}(V)$.

    For the other direction, assume that $df = 0$, and let $u_1, \cdots , u_n$ be a transcendence basis for $\mathbb{C}(V)$ over $\mathbb{C}$. Because $df = 0$, we have that $\partials{}{u_i}(f) = 0$ for all $i$. But since the $u_i$ form a transcendence basis of $\mathbb{C}(V)$ over $\mathbb{C}$, this implies that $f \in \mathbb{C}$.
\end{proof}

\medskip

In the work that follows, we will also need higher arity differential forms, of which we now briefly explain the construction. 

Recall that for any field $F$ and $F$-vector space $U$ with basis $e_1 , \cdots , e_k$ and $m \in \mathbb{N}$, we can form its $m$-th exterior power as the vector space $\bigwedge^m U$ with basis $\{ e_{i_1} \wedge \cdots \wedge e_{i_m} : i_1 < \cdots < i_m\}$. Together, they form the exterior algebra of $U$, i.e the set $\bigwedge U = \bigcup\limits_{m \in \mathbb{N}} \bigwedge^m U$. This a priori only has the structure of an $F$-vector space, but we equip it with the alternating product $\wedge$, which is linear in each coordinate, and satisfies $x \wedge x = 0$ for all $x \in \bigwedge U$. It is an exercise to show that this uniquely defines an alternating product on $\bigwedge U$. 

We can also apply this construction to the dual $U^*$ of $U$. A key observation in this case is that the $m$-th exterior power $\bigwedge^m U^*$ naturally identifies with the $F$-vector space of alternating multilinear function $U^m \rightarrow F$.

We now specialize to the case of the $\mathbb{C}(V)$-vector space $\Omega_V^1$.

\begin{defn}
    For any $m>1$, we define $\Omega^m_V$, the space of \emph{rational $m$-forms} on $V$, as the exterior power $\bigwedge^m \Omega^1_V$. 
\end{defn}

Of particular interest to us will be the space $\Omega_V^{\mathrm{dim}(V)}$. Any element of this space is called a \emph{volume form}. They will be essential tools for classifying algebraic relations between solutions of Lotka-Volterra systems.

By convention, we write $\Omega^0_V = \mathbb{C}(V)$. Recall that $\Omega^1_V$ and $\mathrm{Der}(\mathbb{C}(V)/\mathbb{C})$ are dual spaces. The space $\bigwedge\limits_{i = 1}^m \Omega^1_i$ is the $\mathbb{C}(V)$-vector space of alternating $m$-multilinear maps $\mathrm{Der}(\mathbb{C}(V)/\mathbb{C})^m \rightarrow \mathbb{C}(V)$ (in particular if $m> n$ then $\Omega^m_V = \{ 0 \}$). Concretely, given $1$-forms $\omega_1, \cdots, \omega_m$ and rational sections $s_1,\cdots, s_m$, we have the formula (which we will not use) \[
(\omega_1 \wedge \cdots \wedge \omega_m)(s_1,\ldots,s_m)
=
\sum_{j=1}^m
(-1)^{m+j}\omega_m(v_j)
(\omega_1 \wedge \cdots \wedge \omega_{m-1})
(s_1,\ldots,\widehat{s_j},\ldots,s_m).
\] This sequence of spaces can be made into an exterior algebra, equipped with the alternating product defined above. More precisely, consider the finite dimensional $\mathbb{C}(V)$-vector space \[\Omega_V = \Omega^0_V \oplus \Omega^1_V \oplus \cdots \oplus \Omega^n_V \text{ ,}\]
\noindent the alternating (or wedge) product gives it the structure of a graded $\mathbb{C}(V)$-algebra. Note that for elements of the function field $\Omega_V^0$, the wedge product coincides with the field product. 

Each $\Omega_V^m$ is a vector bundle over $V$, and once again, there is also a local counterpart to this global picture. We can define exterior powers $\bigwedge^m T_p^*V$of the cotangent bundle at a point $p \in V$, which will correspond to alternating multilinear functions $T_p V \rightarrow \mathbb{C}(V)$. If $\omega \in \Omega_V^m$, then $\omega(p) \in  \bigwedge^m T_p^*V$. 

We can extend the exterior derivative canonically to $\Omega_V$. Indeed, one can show that there exists an unique $\mathbb{C}$-linear map extending $d$, satisfying $d \circ d = 0$ and such that for any $k$-form $\omega_1$ and an $l$-form $\omega_2$, we have
\[d(\omega_1 \wedge \omega_2) = d\omega_1 \wedge \omega_2 + (-1)^k \omega_1 \wedge d \omega_2 \text{ .}\]

We therefore have obtained a chain complex
\[0 \rightarrow \Omega^0_V \xrightarrow{d} \Omega^1_V \xrightarrow{d} \cdots \xrightarrow{d} \Omega^n_V \rightarrow 0\]
\noindent of $\mathbb{C}$-vector spaces.

One reason to introduce differential forms, even if one is only interested in vector fields, is that vector fields do not behave well with respect to rational maps, but rational forms do. More precisely, let $\phi : V \rightarrow W$ be a rational map between varieties over $\mathbb{C}$, and $\omega$ a rational $m$-form on $W$. We have a tangent map $T \phi : TV \rightarrow TW$. Fix some $x \in V$, we can define

\begin{align*}
\phi^*_x \omega : (T_x V)^m & \rightarrow \mathbb{C}(V) \\
(v_1 ,\cdots , v_m) & \rightarrow \omega(\phi(x)) (T \phi(v_1) , \cdots , T \phi(v_m)) 
\end{align*}

Because this is an alternating $m$-multilinear map on $T_x V$, we can use it to pullback differential forms:

\begin{defn}
    Let $\phi : V \rightarrow W$ be a rational map between two varieties defined over $\mathbb{C}$. Consider a rational $m$-form over $W$, well-defined on a dense open subset of $W$. We define its \emph{pullback} $\phi^* \omega$ as follows. For all $x \in V$ and rational vector fields $s_1, \cdots , s_m$ defined on $V$, consider, wherever it is well-defined
    \[\phi^* \omega (s_1 ,\cdots , s_m)(x) = \phi^*_x \omega (s_1(x), \cdots , s_m(x)) \text{ .}\]

    \noindent This defines an alternating, $m$-multilinear map
    \begin{align*}
        \phi^* \omega : \mathrm{Der}(\mathbb{C}(V)/\mathbb{C}) & \rightarrow \mathbb{C}(V) \\
        (s_1, \cdots , s_m) & \rightarrow \phi^*\omega(s_1, \cdots , s_m)
    \end{align*}
    i.e. an element of $\Omega_V^m$.
\end{defn}

So far, the structures and functions we have defined depend only on the variety $V$. We are now going to give some structures depending on the choice of a specific rational vector field $s$ on $V$, or equivalently a derivation $D \in \mathrm{Der}(\mathbb{C}(V)/\mathbb{C})$:
\begin{itemize}
    \item the \emph{interior product} $i_s : \Omega^m_V \rightarrow \Omega^{m-1}_V$ defined by $i_s(\omega)(D_1, \cdots , D_{m-1}) = \omega(D, D_1, \cdots , D_{m-1})$. 
    \item the \emph{Lie derivative} $\mathcal{L}_s : \Omega^m_V \rightarrow \Omega^m_V$, which is defined by Cartan's \emph{magic} formula $\mathcal{L}_s = i_D \circ d + d \circ i_D$. Note in particular that for $0$-forms (i.e. elements of the function field), we have that $\mathcal{L}_s = D$.
\end{itemize}

These operators are connected by many formulas, here are a few examples:

\begin{fact}\label{fact: formulas}
    Let $\omega_1 \in \Omega_V^m, \omega_2 \in \Omega_V$, we have:
\begin{itemize}
    \item $\mathcal{L}_s(\omega_1 \wedge \omega_2) = \mathcal{L}_s(\omega_1)\wedge \omega_2+\omega_1 \wedge \mathcal{L}_s(\omega_2)$
    \item $\mathcal{L}_s(df) = d(\mathcal{L}_s(f))$
    \item $i_s(\omega_1 \wedge \omega_2) = i_s(\omega_1) \wedge \omega_2 + (-1)^m \omega_1 \wedge i_s(\omega_2)$
\end{itemize}
If $\omega_1 \in \Omega_V^0$, i.e. $\omega_1 = f \in \mathbb{C}(V)$, then $\omega_1 \wedge \omega_2 = f \omega_2$.
\end{fact}

As an illustration of how these notions all fit together, we provide a simple example computation that will prove useful later:

\begin{prop}\label{prop: int-der-lie-der}
    Let $(V,s)$ be a $D$-variety over some $K < \mathbb{C}$, and pick some $f \in \mathbb{C}(V)$. Then $i_s(df) = \mathcal{L}_s(f)$.
\end{prop}

\begin{proof}
    By definition, we have that $i_s(df) = df(s)$. Recall that $df$ is an element of the dual of $\mathrm{Der}(\mathbb{C}(V)/\mathbb{C})$, obtained by sending $D$ to $D(f)$. The derivation associated to the section $s$ is the Lie derivative $\mathcal{L}_s$, and thus $df(s) = \mathcal{L}_s(f)$. 
\end{proof}

To help our reader develop some intuition, let us explain how the Lie derivative interacts with vector fields. Given a rational vector field $X$ on $V$, the Lie derivative $\mathcal{L}_s(g)$, for some $g \in F$, measures the variation of $g$ along the flow of $s$. This can be made precise using the local existence of integral curves of a vector field, but we give a more algebraic interpretation. Consider a rational vector field on affine space $\mathbb{A}^n$ given by
\[
\begin{cases}
    X_1' = f_1(X_1, \cdots , X_n) \\
    \vdots \\
    X_n' = f_n(X_1, \cdots , X_n)
\end{cases}
\]
\noindent where the $f_i$ are rational functions with complex coefficients. Given any solution $y = (y_1,\cdots,y_n)$ to the system. We can think of $y_i$ as following the flow of the vector field. The Lie derivative of some $g \in \mathbb{C}(y_1, \cdots , y_n)$ is given by
\[\mathcal{L}_s(g) = \sum\limits_{i=1}^n \partials{g}{y_i}f_i\]
\noindent and it is easy to compute that $\mathcal{L}_s(g(y_1, \cdots , y_n)) = D(g(y_1, \cdots , y_n))$, where $D$ is the derivation associated to $s$. Note that $dg$ is the 1-form taking a derivation $D_0$ to $D_0(g)$, and we therefore recover Cartan's formula for the base case of functions.

Given that the Lie derivative measures the variation of a form along the flow of a vector field, the following definition is natural:

\begin{defn}
    Let $s$ be a rational vector field on $V$. We say a differential form $\omega$ is \emph{invariant} if $\mathcal{L}_{s}(\omega) = 0$. The space of invariant $m$-forms is denoted $ \Omega^m_{V,s}$.
\end{defn}

We will give an invariant $1$-form associated to the Lotka-Volterra system:

\[
\begin{cases}
    X'= a_0XY+b_0X \\
    Y' = c_0XY + d_0Y
\end{cases}
\]
This system corresponds to a $D$-variety on $\mathbb{A}^2$ with the section given by
\[s: (x,y) \rightarrow ((x,y),(a_0xy+b_0x,c_0xy+d_0y))\]

\noindent and we have:

\begin{prop}\label{prop: invariant form of LV}
    The $1$-form $\omega = a_0dY - c_0dX + \frac{b_0}{Y}dY - \frac{d_0}{X}dX$ is an invariant form for $\mathcal{L}_s$.
\end{prop}

\begin{proof}
    The essential observation is that
    \begin{align*}
        \mathcal{L}_s(dX) & = d \mathcal{L}_s(X) \text{ by Fact \ref{fact: formulas}}\\
        & = d(b_0 X + a_0 XY) \\
        & = b_0dX + a_0 d(XY) \\
        & = b_0 dX+ a_0YdX + a_0XdY
    \end{align*}
    \noindent and a similar formula for $dY$. We can also compute $\mathcal{L}_s(\frac{1}{X})$ and $\mathcal{L}_s(\frac{1}{Y})$. From there, an elementary computation, using Fact \ref{fact: formulas}, gives the result.
\end{proof}

The form $\omega$ will appear again in the proof of Lemma \ref{lem: complex-inv-curves}.

We recall the following, which is \cite[Lemma 5.6]{freitag2023equations}:

\begin{fact} \label{volume form of dimension 1}
Let \( V \) be an $n$-dimensional affine algebraic variety over \( \mathbb{C} \) and $ D \in \mathrm{Der}(\mathbb{C}(V)/\mathbb{C}) $ be a derivation. If the constant field of $(\mathbb{C}(V),D)$ is equal to $\mathbb{C}$, then the space $\Omega^n_{V,D}$ of invariant volume forms is a one-dimensional complex vector space. 
\end{fact}

A model theorist may find the following useful for intuition: the assumption that the constant field of $(\mathbb{C}(V),D)$ equals $\mathbb{C}$ is equivalent to weak orthogonality to the constants of the generic type of $(V,s)$, where $s$ is the section corresponding to $D$.

We will need the following fact, which is \cite[Proposition 4]{rosenlicht1976liouville}:

\begin{fact}\label{Liouville-Rosenlicht}
    Let $V$ be an irreducible affine variety over $\mathbb{C}$ and consider some $u,v_1, \cdots , v_n \in \mathbb{C}(V)$. If $du + \sum\limits_{i=1}^n c_i\frac{dv_i}{v_i} = 0$ for some $\mathbb{Q}$-linearly independent set $ c_1, \cdots , c_n \in \mathbb{C}$, then $du=0 $ and $dv_i = 0$ for all $1 \leq i \leq n$.
\end{fact}

We note that even without the $\mathbb{Q}$-linear independence, we obtain a similar result that $du$ must be zero, a fact that will be used in the proofs of Theorem \ref{Reproof of Brestovski's Theorem} and Theorem \ref{theo: ortho-desin}:

\begin{cor}\label{cor: liou-ros-not-ind}
    Let $V$ be an irreducible affine variety over $\mathbb{C}$ and consider some $u,v_1, \cdots , v_n \in \mathbb{C}(V)$. Assume that $du + \sum\limits_{i=1}^n c_i\frac{dv_i}{v_i} = 0$ for some $c_1, \cdots , c_n \in \mathbb{C}$. Then $du = 0$.
\end{cor}

\begin{proof}
    If all $c_i$ are zero, the result is immediate, so assume they are not, we can pick $e_1, \cdots , e_m$ to be a maximal $\mathbb{Q}$-linearly independent subset of $c_1, \cdots, c_n$. For any $c_i$, we write $c_i = \sum\limits_{j=1}^m \lambda_{i,j} e_j$ for some $\lambda_{i,j} \in \mathbb{Q}$. Let $N \in \mathbb{Z} \setminus \{ 0 \}$ be such that $\gamma_{i,j} = N \lambda_{i,j}$ is an integer for all $i,j$. 

    For all $i$, we have $Nc_i\frac{d v_i}{v_i} = (\sum\limits_{j=1}^{m} \gamma_{i,j}e_j) \frac{dv_i}{v_i}  = \sum\limits_{j=1}^m e_j \frac{d v_i^{\gamma_{i,j}}}{v_i^{\gamma_{i,j}}}$. Grouping $c_i$ terms together in our first equation, we obtain
    \[N du + \sum\limits_{j=1}^m e_j \sum\limits_{i=1}^n \frac{d v_i^{\gamma_{i,j}}}{v_i^{\gamma_{i,j}}} = 0\]
    \noindent and from there that
    \[du + \sum\limits_{j=1}^m \frac{e_j}{N} \frac{d w_j}{w_j} = 0\]

    \noindent where $w_j = \prod\limits_{i=1}^n v_i^{\gamma_{i,j}}$. This yields the result by Fact \ref{Liouville-Rosenlicht}.
\end{proof}

\subsection{A sufficient criterion for strong minimality}
In this subsection, we present the general method we will use to prove strong minimality of the Lotka-Volterra systems (both classical and 2d). Note that the method that we will end up using, Corollary \ref{cor: criterion for sm}, could be deduced more or less directly from Jaoui's \cite[Theorem D]{jaoui2023density}, but we chose an alternative path to highlight some potentially useful independent result. 

The main idea is to connect strong minimality and proper almost fibrations. The notion of proper fibrations was introduced in \cite{moosa2014some} and we now recall its definition. 

\begin{defn}
    A type $p$ over some algebraically closed differential field $K$ \emph{admits a proper fibration} if there is some realization $a$ of $p$ and some $b$ such that $b \in K \langle a \rangle \setminus K$ and $a \not\in K\langle b \rangle^{\mathrm{alg}}$. If instead we have $b \in K \langle a \rangle^{\mathrm{alg}} \setminus K$, we say that $p$ admits a proper \emph{almost} fibration.
\end{defn}

Observe that if a type has a proper almost fibration, then it is not minimal (and thus not strongly minimal). We are now going to explain why the converse is true for generic types of planar rational vector fields: if they have no proper fibrations, then they must be strongly minimal (note that this is false in general, examples are given in \cite[Section 4]{freitag2023bounding}). 

We use the following criteria for minimality:

\begin{thm}\label{theo: min-crit}
    Let $K$ be an algebraically closed field of constants, some $f_1, \cdots f_k \in K(X_1, \cdots , X_k)$ and $p \in S(K)$ the generic type of the system:
    \[\begin{cases}
    X_1' = f_1(X_1, \cdots , X_k) \\
    \vdots \\
    X_k' = f_k(X_1, \cdots , X_k)
    \end{cases}\]
    
    \noindent The type $p$ is minimal if and only if it has no proper almost fibrations. 
\end{thm}

\begin{proof}
    The left to right direction is immediate. For the other direction, assume that $p$ has no proper almost fibrations. We may assume that $k > 1$ (if $k = 1$ then $p$ is minimal). By \cite[Theorem 6.8]{jaoui2022abelian}, either $p$ is minimal, or there is a finite-to-one $K$-definable map $f : p \rightarrow q$ where $q$ is the generic type of a logarithmic-differential equation on a simple abelian $A$ variety of dimension strictly greater than one (see \cite[Section 6.2]{jaoui2022abelian} for a definition). This implies that the binding group of $q$ is definably isomorphic to a definable subgroup of the $\mathbb{C}$-points of $A$. In particular, if it is infinite, then it cannot be linear. 

    However, by \cite[Lemma 3.1]{eagles2024internality}, the binding group of $q$ must be linear, and thus finite by the previous discussion. This means that any realization of $q$ is algebraic over $\mathbb{C}$, and as $\mathbb{C}$ is algebraically closed, $q$ must be the type of a tuple of constants. But since $k > 1$, it has transcendence degree strictly greater than $1$, and thus must have a proper almost fibration, which is a contradiction. 
\end{proof}

To prove that a type has no proper almost fibrations, we will use the following key result of Jaoui, which is \cite[Theorem B]{jaoui2022generic}:

\begin{fact}[Jaoui]\label{theo: Jaoui B}
    Let $(V,s)$ be a $D$-variety of dimension $2$, defined over constant parameters. If there is a singular point $\overline{x} \in D(\mathbb{C})$ such that $\overline{x}$ is not contained in any complex invariant algebraic curve, then the generic type of $(V,s)$ admits no proper almost fibration. 
\end{fact}

Putting these results together, we obtain:

\begin{cor}\label{cor: criterion for sm}
    Let $f,g \in K[X,Y]$, for some algebraically closed $K < \mathbb{C}$, and consider the definable set
    \[\label{order 2 sys}\tag{$S$}
    \begin{cases}
        X' = f(X,Y) \\
        Y' = g(X,Y) \text{ .}
    \end{cases}
    \]
    \noindent If the associated algebraic vector field has a singular point that is not contained in any invariant curve defined over $\mathbb{C}$, then the generic type of (\ref{order 2 sys}) is strongly minimal.
\end{cor}

\begin{proof}
    We obtain minimality simply by applying Theorem \ref{theo: Jaoui B} and Corollary \ref{cor: criterion for sm}. From there, we get strong minimality of the generic type as for order two definable sets, Morley rank and Lascar rank coincide (see \cite[Theorem 6.1]{freitag2017finiteness}). 
\end{proof}

\subsection{A generalization of a theorem of Brestovski}\label{subsec: brestovski}

In this subsection, we generalise the following theorem of Brestovski, which is \cite[Theorem 2]{brestovski1989algebraic}. Fix some countable algebraically closed field of constants $K < \mathbb{C}$ (we have to pick it countable so that is is small with respect to our ambient model $\mathcal{U}$). 

\begin{fact}[Brestovski]
Let $F,G_1,...,G_n\in \mathbb{C}(X_0,X_1)$ be rational functions and $b_1,...,b_n\in\mathbb{C}$  be $\mathbb{Q}$-linearly independent. Assume that the generic type of the set of solutions to the equation in $K\langle X \rangle$
\[
F(X,X')'=\sum\limits_{i=1}^nb_i\frac{G_i(X,X')'}{G_i(X,X')}
\]
\noindent is minimal. For any two generic solutions $a_1$ and $a_2$ of this equation, if $a_1 \in K\langle a_2 \rangle^{\mathrm{alg}}$, then there is a non-zero $m \in \mathbb{N}$ such that $F(a_1)^m = F(a_2)^m$
\end{fact}

Our theorem will generalize this to classify relations between two equations of this form. It also removes the $\mathbb{Q}$-linear independence assumption between the coefficients. 

However, note that our result is slightly weaker: in his proof, Brestovski uses the intermediate result that there are constants $e,f \in \mathbb{C}$ such that $F(a_1) = eF(a_2) + f$. We could only generalize his work as far as this equation. Our proof more or less follows Brestovski's, but is written in a more modern language.

\begin{thm}\label{Reproof of Brestovski's Theorem} 
Let $F_1, G_1,...,G_n\in K(X_0,X_1)$ and $F_2,H_1,...,H_m\in K(Y_0,Y_1 )$ be rational functions. Let $b_1,...,b_n,c_1,...,c_m\in K$. 

Consider the following equations in $K\langle X \rangle$ and $K \langle Y \rangle$
\begin{equation} \tag{$S_1$} \label{$S_1$}
   F_1(X,X')' = \sum_{i=1}^n b_i\frac{G_i(X,X')'}{G_i(X,X')}
\end{equation}
\begin{equation} \tag{$S_2$} \label{$S_2$}
   F_2(X,X')' = \sum_{i=1}^m c_i\frac{H_i(X,X')'}{H_i(X,X')}.
\end{equation}

and assume that the generic types $p_1$ and $p_2$ of their solution sets are minimal. 

If $a_1$ is a generic solution of $(S_1)$ and $a_2$ is a generic solution of $(S_2)$ such that $a_1\in K\langle a_2\rangle^{\mathrm{alg}}$, then there are $e,f \in \mathbb{C}$ such that $F_1(a_1,a_1') = eF_2(a_2,a_2') +f$.
\end{thm}

\begin{proof}
    By clearing the denominators, we see that $(S_1)$ has the same solution set as an equation of the form $X'' P_1(X,X') + Q_1(X,X')=0$, for some $P_1,Q_1 \in K[X_0,X_1]$. We similarly obtain polynomials $P_2,Q_2 \in K[Y_0,Y_1]$ by clearing the denominators in $(S_2)$. By the discussion following Fact \ref{fact: d-var-capture-types}, we can thus identify the solution sets of $(S_1)$ and $(S_2)$ to the sharp points of rational $D$-varieties $(\mathbb{A}^2,s_1)$ and $(\mathbb{A}^2, s_2)$, respectively, where the section $s_1$ is given by
    \begin{align*}
        s_1 : \mathbb{A}^2 & \rightarrow T \mathbb{A}^2 \\
        (x_0, x_1) & \rightarrow \left( (x_0,x_1), \left( x_1, \frac{-Q_1(x_0,x_1)}{P_1(x_0,x_1)} \right)\right)
    \end{align*}
    and a similar formula for $s_2$. In particular, we see that $\mathcal{L}_{s_1}(X_0) = X_1$ and $\mathcal{L}_{s_1}(X_1) = \frac{-Q_1(X_0,X_1)}{P_1(X_0,X_1)}$. Again, we have the same formulas for $\mathcal{L}_{s_2}$. The intuition here is that the variables $(X_0,X_1)$ are identified with $(X,X')$, and $X''$ is identified with $\frac{-Q_1(X_0,X_1)}{P_1(X_0,X_1)}$. To be rigorous, we will only use $X_0,X_1, Y_0$ and $Y_1$ in our proof, but the reader should keep these identifications in mind. 

    From the equation $(S_1)$, we obtain that the Lie derivative $\mathcal{L}_{s_1}$ satisfies the following equation on $K(X_0,X_1)$
    \begin{equation}
        \tag{$\mathcal{L}S_1$} \mathcal{L}_{s_1}(F_1) = \sum\limits_{i=1}^n b_i \frac{\mathcal{L}_{s_1}(G_i)}{G_i}
    \end{equation}
    \noindent and a similar equation $(\mathcal{L}S_2)$ for $\mathcal{L}_{s_2}$.

    Using these equations, we will now identify invariant volume forms for our rational $D$-varieties $(\mathbb{A}^2,s_1)$ and $(\mathbb{A}^2,s_2)$. Consider the volume forms
    \begin{align*}
        \omega_1 & = \frac{dF_1}{\mathcal{L}_{s_1}(F_1)} \wedge \sum_{i=1}^nb_i\frac{dG_i}{G_i}\\
        \omega_2 & = \frac{dF_2}{\mathcal{L}_{s_2}(F_2)} \wedge \sum_{i=1}^m c_i\frac{dH_i}{H_i}
    \end{align*}

    We start by computing their interior derivatives with respect to $s_1$ and $s_2$:

    \begin{claim}\label{claim: int-der-comp}
        We have the formulas
        \begin{align*}
           i_{s_1}(\omega_1) & =   \sum\limits_{i=1}^n \frac{dG_i}{G_i} - dF_1\\
           i_{s_2}(\omega_2) & = \sum\limits_{i=1}^m \frac{dH_i}{H_i} - dF_2
        \end{align*}
    \end{claim}

    \begin{proof}[Proof of claim]
        We only prove this for $\omega_1$, the computation for $\omega_2$ being identical. Recall that by Proposition \ref{prop: int-der-lie-der}, for any $f \in \mathbb{C}(\mathbb{A}^2) = \mathbb{C}(X_0,X_1)$, we have $i_{s_1}(df) = \mathcal{L}_{s_1}(f)$. Using this and Fact \ref{fact: formulas}, we compute:
        \begin{align*}
            i_{s_1}(\omega_1)
            & = i_{s_1}\left( \frac{dF_1}{\mathcal{L}_{s_1}(F_1)}\right)\wedge \sum\limits_{i=1}^n \frac{dG_i}{G_i}- \frac{dF_1}{\mathcal{L}_{s_1}(f_1)} \wedge\sum\limits_{i=1}^n b_i i_{s_1}\left( \frac{dG_i}{G_i}\right) \\
            & =  \frac{\mathcal{L}_{s_1}(F_1)}{\mathcal{L}_{s_1}(F_1)} \wedge \sum\limits_{i=1}^n \frac{dG_i}{G_i} - \frac{dF_1}{\mathcal{L}_{s_1}(F_1)} \wedge \sum\limits_{i=1}^n b_i  \frac{\mathcal{L}_{s_1}(G_i)}{G_i}  \text{ by Proposition \ref{prop: int-der-lie-der}} \\
            & =  \sum\limits_{i=1}^n \frac{dG_i}{G_i} - \frac{dF_1}{\mathcal{L}_{s_1}(F_1)} \wedge \sum\limits_{i=1}^n b_i  \frac{\mathcal{L}_{s_1}(G_i)}{G_i} \\
            & = \sum\limits_{i=1}^n \frac{dG_i}{G_i} - \frac{dF_1}{\mathcal{L}_{s_1}(F_1)} \wedge \mathcal{L}_{s_1}(F_1)  \text{ by equation $(\mathcal{L}S_1)$} \\
            & = \sum\limits_{i=1}^n \frac{dG_i}{G_i} - dF_1 
        \end{align*}
    \end{proof}

    We can now quickly prove that the forms $\omega_1$ and $\omega_2$ are invariant. Again, we only give the proof for $\omega_1$. First note that Cartan's formula gives us $\mathcal{L}_{s_1}(\omega_1) = (d \circ i_{s_1} + i_{s_1} \circ d)(\omega_1)$. Because $\omega_1$ is a volume form, the second summand cancels and we get $\mathcal{L}_{s_1}(\omega_1) = (d \circ i_{s_1})(\omega_1)$. By the claim, we get
    \begin{align*}
        \mathcal{L}_{s_1}(\omega_1) & = (d \circ i_{s_1})(\omega_1) \\
        & = d\left(\sum\limits_{i=1}^n \frac{dG_i}{G_i} - dF_1\right) \\
        & = \sum\limits_{i=1}^n d\left(\frac{dG_i}{G_i}\right) \text{ as $d \circ d = 0$} \\
        & = \sum\limits_{i=1}^n d(dG_i)\wedge \frac{1}{G_i} - dG_i \wedge d\left(\frac{1}{G_i}\right) \\
        & = - \sum\limits_{i=1}^n dG_i \wedge d\left(\frac{1}{G_i}\right) \\
        & = \sum\limits_{i=1}^n \frac{dG_i \wedge dG_i}{G_i^2} \text{ as $d\left( \frac{1}{G_i} \right) = -\frac{dG_i}{G_i^2}$}\\
        & = 0 \text{ as $dG_i \wedge dG_i = 0$ .}
    \end{align*}

    Now assume that there are $a_1 \models p_1$ and $a_2 \models p_2$ such that $a_2 \in K\langle a_2 \rangle$. Then $(a_1,a_1')$ and $(a_2,a_2')$ are generic sharp points of $(\mathbb{A}^2,s_1)$ and $(\mathbb{A}^2,s_2)$, respectively. Since $p_1$ and $p_2$ are minimal, we have that $(a_1,a_1')$ and $(a_2,a_2')$ are interalgebraic, i.e. $a_2,a_2' \in K(a_1,a_1')^{\mathrm{alg}}$ and vice-versa. Consider now the Zariski locus $W$ of $(a_1,a_1',a_2,a_2')$ in $\mathbb{A}^4$. Because of interalgebraicity, this is a subvariety of $\mathbb{A}^4$ of dimension $2$. Moreover, the map $s_1 \times s_2 : \mathbb{A}^4 \rightarrow T\mathbb{A}^4$ restricts to a section of the tangent bundle of $W$, and we thus obtain the $D$-variety $(W,s_1 \times s_2)$. Note that we also have the sections $s_1 \times 0$ and $0 \times s_2$ of $TW$, which we will denote $s_1$ and $s_2$ in the rest of the proof.

    Let $\pi_1$ and $\pi_2$ be the coordinate projections on the first and last two coordinates, respectively. We obtain pullback forms $\pi_1^* \omega_1$ and $\pi_2^* \omega_2$ on $W$, which we still denote, by a slight abuse of notation, $\omega_1$ and $\omega_2$. They are invariant on $(W,s_1 \times s_2)$ because $\omega_1$ and $\omega_2$ are invariant on $(\mathbb{A}^2,s_1)$ and $(\mathbb{A}^2,s_2)$, respectively. Because $p_1$ and $p_2$ are minimal and of order strictly greater than one, by Fact \ref{fact: mini-implies-ortho} they are orthogonal to the constants. Therefore $\tp(a_1,a_2/K)$ is orthogonal to the constants. By Fact \ref{volume form of dimension 1}, there is $e \in \mathbb{C}$ such that $\omega_1 = e \omega_2$.

    Consider the differential form on $W$ given by
    \[\omega = d(eF_2-F_1) + \sum\limits_{i=1}^n b_i \frac{d G_i}{G_i}-\sum\limits_{i=1}^nec_i \frac{dH_i}{H_i}\]
    \noindent we now show that $\omega = 0$. Since $\omega_1 - e \omega_2 = 0$, we also have $i_{s_1 \times s_2}(\omega_1 - e \omega_2) = 0$. Therefore
    \begin{align*}
        0 & = i_{s_1 \times s_2}(\omega_1-e \omega_2) \\
        & = i_{s_1 \times s_2}(\omega_1) -e i_{s_1 \times s_2}(\omega_2) \\
        & = i_{s_1}(\omega_1) -
        ei_{s_2}(\omega_2) \\
        & = \sum\limits_{i=1}^n \frac{dG_i}{G_i} - dF_1 -e \left( \sum\limits_{i=1}^m \frac{dH_i}{H_i} - dF_2\right) \text{ by Claim \ref{claim: int-der-comp}}\\
        & = \omega
    \end{align*}
    We can now apply Corollary \ref{cor: liou-ros-not-ind} to conclude that $d(eF_2-F_1) =0$, and thus by Proposition \ref{prop: df-zero-implies-cst}, that $eF_2 - F_1 = f$ for some $f \in \mathbb{C}$.
\end{proof}

\section{Lotka-Volterra systems}\label{sec: LV}

In this section, we will assume that $a,b,c,d \neq 0$. We use the following notation to denote the Lotka-Volterra system of equations
   \begin{equation}\notag
LV_{a,b,c,d}:=\begin{dcases}
    X'= X(aY + b) \\
    Y'= Y(cX + d)
\end{dcases}
\end{equation}

Our goal is to entirely classify the differential-algebraic relations between generic solutions of Lotka-Volterra systems. 

Let us give some useful transformations between different Lotka-Volterra systems. First, we have the rescaling transformation
\begin{align*}
    LV_{a,b,c,d} & \rightarrow LV_{1,b,1,d} \\
    (x,y) & \rightarrow (cx,ay)
\end{align*}
\noindent which induces a bijection between solutions of these two systems.

We also have, because of the symmetry of the system, the swapping transformation
\begin{align*}
    LV_{a,b,c,d} & \rightarrow LV_{c,d,a,b} \\
    (x,y) & \rightarrow (y,x)
\end{align*}
\noindent which again is a bijection between solutions of the two systems.

As we will see in Subsection \ref{subsec: ortho}, these are the only possible relations.

\subsection{Complex invariant curves and strong minimality}\label{subsec: complex inv curve}

In this subsection, we will take the first key step of identifying the complex invariant curves. Bearing in mind Corollary \ref{cor: criterion for sm}, this will immediately give strong minimality.

In the $b \neq d$ case of the classical system, the complex invariant curves were already determined by Duan and Nagloo in \cite[Corollary 3.4]{duan2025algebraic}. We include a different proof that also deals with the $b=d$ case, and is a good illustration of the use of rational forms in this context.

Recall that the vector field associated to the Lotka-Volterra system $LV_{a,b,c,d}$ is
\begin{align*}
    s: \mathbb{A}^2 & \rightarrow T\mathbb{A}^{2} \\
    (x,y) & \rightarrow \left( (x,y) , (bx + axy, dy + cxy)\right)
\end{align*}

\begin{lem}\label{lem: complex-inv-curves}
    If $b \neq d$, the complex irreducible invariant curves of the vector field $s$ associated to $LV_{a,b,c,d}$ are given by $X= 0$ and $Y =0$. If $b = d$, there is also $cX - aY = 0$. 
\end{lem}

\begin{proof}
    The fact that these are invariant curves is a quick computation, left to the reader.

    By using the first of the previously mentioned transformations, we can reduce to the case of $LV_{1,b,1,d}$ (we will just have to be mindful of that change for the case $b=d$). In the rest of the proof, we will denote $b$ and $d$ by $b_0$ and $d_0$, to avoid confusion with the exterior derivative $d$.

    We will use the the rational differential $1$-form $\omega = dY-dX+ \frac{b_0}{Y}dY - \frac{d_0}{X}dX$, defined outside of $XY = 0$. Note that $s$ is a section of the tangent bundle of $\mathbb{A}^2$, and thus identified with an element of $\mathrm{Der}(\mathbb{C}(X,Y)/\mathbb{C})$. More precisely, it is identified with the derivation $f \rightarrow \mathcal{L}_s(f)$. Also recall that $\omega \in \Omega_{\mathbb{A}^2}^1$, which is defined to be the dual of $\mathrm{Der}(\mathbb{C}(X,Y)/\mathbb{C})$ as a $\mathbb{C}(X,Y)$ vector space. We can thus apply $\omega$ to $s$ to obtain an element of $\mathbb{C}(X,Y)$. We will now show that $\omega(s) = 0$.
    First note that
    \begin{align*}
        dX(s) & = s(X) \\
        & = \mathcal{L}_s(X) \\
        & = b_0 X + XY
    \end{align*}
    and similarly $dY(s) = d_0 Y + XY$. From this we compute that
    
    \begin{align*}
        \omega(s) & = dY(s) -dX(s) + \frac{b_0}{Y}dY(s) - \frac{d_0}{X}dX(s) \\
        & = d_0 Y + XY - b_0X-XY + \frac{b_0}{Y}(d_0 Y + XY) - \frac{d_0}{X}(b_0 X + XY) \\
        & = 0
    \end{align*}

    Now fix a complex invariant irreducible curve $C$, different from $X = 0$ and $Y = 0$. We have an inclusion map $C \rightarrow \mathbb{A}^2$, and we can consider the pullback $\omega^*$ along this map, it is an element of $\Omega_C^1$. We will now argue that $\omega^* = 0$. 

    Let $r$ be a rational section of $TC$, which we can identify with an element of $\mathrm{Der}(\mathbb{C}(C)/\mathbb{C})$. Note that $\mathrm{Der}(\mathbb{C}(C)/\mathbb{C})$ is a one-dimensional $\mathbb{C}(C)$-vector space, and thus there is $\lambda \in \mathbb{C}(C)$ such that $r = \lambda s$. From this we compute, for any $p \in C$ at which $r$ is well-defined, that
    \begin{align*}
        \omega^*(r)(p) & = \omega(p)(r(p)) \\
        & = \omega(p)(\lambda(p) s(p)) \\
        & = \lambda(p) \omega(s)(p) \\
        & = 0
    \end{align*}
    \noindent which shows that $\omega^* = 0$.

    From $\omega^* = 0$, we conclude that in $\Omega^1(\mathbb{C}(C)/\mathbb{C})$, we have
    \[d(Y-X) +  \frac{b_0}{Y}dY -  \frac{d_0}{X}dX = 0\text{ .}\]

    By Corollary \ref{cor: liou-ros-not-ind}, we see that $d(Y-X) = 0$. In other words $Y-X = \lambda$ for some $ \lambda \in \mathbb{C}$, implying that $C$ is contained in the curve $Y-X - \lambda$, and by irreducibility that $C$ equals that curve. 

    Therefore $Y-X -\lambda$ is an invariant curve, which by definition means that $Y-X-\lambda$ divides $\mathcal{L}_s(Y-X-\lambda)$, as $\mathcal{L}_s$ is the derivation associated with $s$. A quick computation yields $\mathcal{L}_s(Y-X-\lambda) = d_0Y-b_0X$, and from there we can conclude that we must have $b_0 = d_0$ and $\lambda = 0$.

    Remembering that we used the transformations $x \rightarrow cx$ and $y \rightarrow ay$, we obtain the lemma.
    
\end{proof}

\begin{rem}
    The reader familiar with differential forms and the Lotka-Volterra system may have noticed a quicker way to prove that $\omega(s) = 0$. Indeed, note that there is a first integral for the vector field $s$ (with $a_0 = c_0 = 1$) given by
    \[f = Y-X+b_0 \log(Y) - d_0 \log(X)\]
    meaning that $\mathcal{L}_s(f) = 0$. Moreover $df = \omega$. Then one computes that $\omega(s) = df (s) = \mathcal{L}_s(f) = 0$. We chose to give a purely algebraic proof to avoid introducing even more material.
\end{rem}

For the 2d Lotka-Volterra system
\[
\begin{cases}\tag{$LV_{a,b,c,d}^{2d}$}
    X' = X(aY +b) \\
    Y' = Y(cX + dY)
\end{cases}
\]
\noindent with $a,b,c,d \in \mathbb{C} \setminus \{ 0 \}$, we have the following, which is \cite[Corollary 3.8]{duan2025algebraic}:

\begin{fact}
    The complex invariant curves of the 2d Lotka-Volterra system are exactly $X = 0$ and $Y = 0$.
\end{fact}

\begin{defn}
    We will call a solution of $LV_{a,b,c,d}$ or $LV_{a,b,c,d}^{2d}$ \emph{degenerate} if it lies in any of these complex invariant curves.
\end{defn}

Note that this allows us, in particular, to identify the generic types of $LV_{a,b,c,d}$ and $LV_{a,b,c,d}^{2d}$:

\begin{cor}
    The generic types of $LV_{a,b,c,d}$ and $LV_{a,b,c,d}^{2d}$ are the unique types of non-degenerate, non-constant solutions.
\end{cor}
\begin{proof}
    Seeing the set defined by $LV_{a,b,c,d}$ as a $D$-variety on $\mathbb{A}^2$, its generic type is simply the type of a sharp point not in any proper $D$-subvariety, or equivalently, non-constant and not in any invariant curve. The same goes for $LV_{a,b,c,d}^{2d}$.
\end{proof}

From these classifications of complex invariant algebraic curves and Corollary \ref{cor: criterion for sm}, we deduce:

\begin{thm}\label{theo: LV-sm}
    Let $a,b,c,d \in \mathbb{C} \setminus \{ 0 \}$. The definable set of non-degenerate solutions of the Lotka-Volterra system
    \[\begin{cases}
        X' = aXY + bX \\
        Y' = cXY + dY \\
        X \neq 0 \\
        Y \neq 0
    \end{cases}\]
    is strongly minimal if and only if $b \neq d$.
\end{thm}

\begin{proof}
    That it is not strongly minimal if $b = d$ is well known, see \cite[subsubsection 3.1.1]{duan2025algebraic} for example. For the other direction, the singular point $(\frac{-d}{c}, \frac{-b}{a})$ does not lie on the invariant curves $X= 0$ and $Y=0$.
\end{proof}

\begin{rem}
    It is clear why the proof does not work if $b=d$: in that case, the singular point $(\frac{-d}{c},\frac{-b}{a})$ \emph{does} lie on the third invariant curve $cX-aY =0$.
\end{rem}

Duan and Nagloo already proved strong minimality if $\frac{b}{d} \not\in \mathbb{Q}$ in \cite[Theorem A]{duan2025algebraic}. Our new contribution is the $\frac{b}{d} \in \mathbb{Q}$, but $b \neq d$ case.

As for the $2d$-system, Duan and Nagloo, in \cite[Theorem B]{duan2025algebraic}, proved strong minimality of the definable set of non-degenerate solutions. We do point out how we can now quickly recover their result:

\begin{thm}[Duan-Nagloo]
The definable set of non-degenerate solutions of the 2d Lotka-Volterra system
    \[
    \begin{cases}
    X' = X(aY + b) \\
    Y' = Y(cX + dY)\\
    X \neq 0 \\
    Y \neq 0
    \end{cases}
    \]
\noindent with $a,b,c,d \in \mathbb{C} \setminus \{ 0 \}$ is strongly minimal.
\end{thm}

\begin{proof}
    Here, we just have to observe that the only singular point is $(\frac{cd}{ca}, \frac{-b}{a})$, which is not in the complex invariant curves $X=0$ and $Y = 0$.
\end{proof}

\subsection{Orthogonality and disintegration}\label{subsec: ortho}

Our goal in this section is to entirely classify the algebraic relations between solutions of Lotka-Volterra systems, as long as $b \neq d$. 

We start by proving that the only relation between solutions of a fixed Lotka-Volterra system is equality (total disintegration) and that the relations between different systems reduce to the two types given in subsection \ref{subsec: complex inv curve}. As a reminder, these are the rescaling transformation
\begin{align*}
    LV_{a,b,c,d} & \rightarrow LV_{1,b,1,d} \\
    (x,y) & \rightarrow (cx,ay)
\end{align*}
and the swapping transformation
\begin{align*}
    LV_{a,b,c,d} & \rightarrow LV_{c,d,a,b} \\
    (x,y) & \rightarrow (y,x)
\end{align*}

In the next few results, we only work with systems of the form $LV_{1,b,1,d}$ for simplicity. For these systems, the next theorem shows that the swapping transformation is the only possible algebraic relation (besides equality). 

\begin{thm}\label{theo: ortho-desin}
    Let $LV_{1,b_1,1,d_1}$ and $LV_{1,b_2,1,d_2}$ be Lotka-Volterra systems with $\frac{b_i}{d_i} \neq 1$ for $i = 1,2$. Let $(x,y)$ and $(u,v)$ be generic solutions of $LV_{1,b_1,1,d_1}$ and $LV_{1,b_2,1,d_2}$, respectively. If $(x,y)$ and $(u,v)$ are not independent over $\mathbb{C}$, then either:
    \begin{itemize}
        \item $(b_1,d_1) = (b_2,d_2)$ and $(x,  y) = (u,v)$,
        \item $(b_1,d_1) = (d_2,b_2)$ and $(x,y) = (v,u)$.
    \end{itemize}
\end{thm}

\begin{proof}
    Let $(x,y)$ and $(u,v)$ be generic solutions of $LV_{1,b_1,1,d_1}$ and $LV_{1,b_2,1,d_2}$ respectively. First note that because the generic types of these equations are orthogonal to the constants by Fact \ref{fact: mini-implies-ortho}, they are also generic solutions over $\mathbb{C}$, meaning that $\mathrm{trdeg}(x,y/\mathbb{C})= \mathrm{trdeg}(u,v/\mathbb{C}) = 2$. Assume that $(x,y)$ and $(u,v)$ are not independent over $\mathbb{C}$, then there exists a countable subfield $K < \mathbb{C}$ such that they are not independent over $K$. By strong minimality this implies that they are interalgebraic over $K$.
    
    Consider $z = x-y$ and $w = u-v$.
    It is easy to compute that
    \[z' = b_1x - d_1 y = b_1\frac{y'}{y} - d_1\frac{x'}{x}\]
    and that
    \[
    \begin{cases}
        z'- d_1 z = (b_1-d_1)x \\
        z'- b_1z = (b_1-d_1)y
    \end{cases}
    \]
    from which we can deduce, as $b_1\neq d_1 \neq 1$, that $z$ is interdefinable with $(x,y)$ over $\mathbb{Q}(b_1,d_1)$ and satisfies the differential equation:
    \[Z'=b_1\frac{(Z'-b_1Z)'}{Z'-b_1Z}- d_1\frac{(Z'-d_1 Z)'}{Z'- d_1 Z} \text{ .}\]
    Similarly, we see that $w$ is interdefinable with $(u,v)$ over $\mathbb{Q}(b_2,d_2)$ and satisfies the equation:
    \[W'=b_2\frac{(W'-b_2 W)'}{W'-b_2W}- d_2\frac{(W'-d_2 W)'}{W'- d_2 W} \text{ .}\]
    By interdefinability, we get that $z$ and $w$ are interalgebraic over $K$. The two equations they satisfy are of Brestovski's form, therefore we can apply Theorem \ref{Reproof of Brestovski's Theorem} (interdefinability yields that $z$ and $u$ are generic solutions). Thus there are constants $e,f$ such that we have $z = ew+f$ or in terms of $x,y,u,v$:

    \begin{align}\label{O1}
    x-y = e(u-v) +f \text{ .}
    \end{align}

    Note that this implies $e \neq 0$, as  otherwise $x-y = f$, which is impossible as $(x,y)$ is a generic solution of $LV_{1,b_1,1,d_1}$. Differentiating (\ref{O1}) we have 
    \begin{align}\label{O2}
        b_1x- d_1 y & = eb_2u - ed_2 v 
    \end{align}
    and the operations $d_1 \eqref{O1} - \eqref{O2}$ and $b_1\eqref{O1} - \eqref{O2}$ give us a transformation between solutions, that is:
    \begin{equation}
    \begin{dcases}
        x  = Au + Bv + d_1G \\
        y  = Cu + Dv + b_1G\\
    \end{dcases}
    \end{equation} where 
    $A = \frac{b_2-d_1}{b_1-d_1}e$, $B = \frac{d_1-d_2}{b_1-d_1}e $, $ C = \frac{b_2-b_1}{b_1-d_1}e $, $ D = \frac{b_1-d_2}{b_1-d_1}e$ and $G = \frac{-f}{b_1-d_1}$.

    Using these identities and the Lotka-Volterra equations, we can compute $x'$ two ways and obtain the following polynomial equality in $u,v$:
    \begin{align*}
        ACu^2 & + BDv^2 + (AD+CB-A-B)uv\\
        & +(b_1AG+d_1CG +b_1A-b_2A)u \\
        & + (d_1GD+b_1GB+b_1B-d_2B)v\\
        & + b_1d_1G^2 + b_1d_1G = 0 \text{ .}
    \end{align*}
    If the coefficients are not all zero, this implies that the transcendence degree of $(u,v)$ over the field generated by $b_1,d_1,b_2,d_2, e$ and $f$ is less than one. By strong minimality of the system $LV_{1,b_2,1,d_2}$, it must be zero, which is impossible as the generic type of $LV_{1,b_2,1,d_2}$ is orthogonal to the constants. 

    Therefore all coefficients of this polynomial must be zero. The coefficients of $u^2$ and $v^2$ give us $AC = 0$ and $BD = 0$. Going back to the definitions of $A,B,C$ and $D$, this means that $b_2=d_1$ or $b_2 = b_1$ and also $d_1 = d_2$ or $b_1 = d_2$. Using that $b_1 \neq d_1$ and $b_2 \neq d_2$, we reduce to the two cases:
    \[\begin{cases}
        b_1 = d_2 \text{ and } b_2 = d_1 \\
        b_1 = b_2 \text{ and } d_1 = d_2 \text{ .}
    \end{cases}\]
    Note that $(x,y) \rightarrow (y,x)$ gives a definable bijection between $LV_{1,b_1,1,d_1}$ and $LV_{1,d_1,1,b_1}$. Therefore, we can reduce the first case to the second, and we now assume that $b_1 = b_2$ and $d_1 = d_2$, i.e. $C = B = 0$. Again going back to the definitions, this immediately implies that $A = D$.

    Examining the $uv$ terms, this leads to $A^2 - A = 0$, so either $A = 0$ or $A = 1$. The first case is impossible, as it would lead to $b_2 = d_2$. Therefore $A = 1$ and $D = 1$.

    Finally, examining the $v$ term, we see that $d_1 G = 0$, and as $d_1 \neq 0$ we get $G = 0$, so $f = 0$.

    Keeping in mind that we may have swapped $x$ and $y$, we obtain the result.

\end{proof}

In terms of types, we have the following consequence:

\begin{cor}
    Let $LV_{1,b_1,1,d_1}$ and $LV_{1,b_2,1,d_2}$ be Lotka-Volterra systems with $\frac{b_i}{d_i} \neq 1$ for $i = 1,2$. Then
    \begin{enumerate}
        \item their generic types are totally disintegrated,
        \item if $\{ b_1, d_1\} \neq \{ b_2, d_2\}$, their generic types are orthogonal.
    \end{enumerate}
\end{cor}

\begin{proof}
    Let $p_1$ and $p_2$ be their generic types. In the previous theorem, we have obtained that if $(x,y)$ and $(u,v)$ are realizations of $p_1$ and $p_2$, not independent over $K = \mathbb{Q}(b_1,d_1,b_2,d_2)^{\mathrm{alg}}$, then either $(b_1,d_1) = (b_2,d_2)$ and $(x,y) = (u,v)$, or $(b_1,d_1) = (d_2,b_2)$ and $(x,y) = (v,u)$. In particular, $p_1$ and $p_2$ are totally disintegrated: the only relation between solutions of $p_1$ is equality, and the same is true for $p_2$. Moreover, if $\{ b_1, d_1\} \neq \{ b_2,d_2\}$ and $(x,y)$ and $(u,v)$ realize $p_1$ and $p_2$, respectively, they must be independent. Therefore $p_1$ and $p_2$ are weakly orthogonal. Thus they must be orthogonal by Fact \ref{weak non-orthogonality}, as they are disintegrated. 
\end{proof}

We can now prove our main result, which completely classifies algebraic relations between generic solutions of \emph{potentially different} Lotka-Volterra systems:

\begin{thm}\label{theo: main-theo}
    Let $LV_{a_i,b_i,c_i,d_i}, 1 \leq i \leq n$ be Lotka-Volterra systems 
    \[
    \begin{cases}
        X' = a_i XY + b_iX  \\
        Y' = c_i XY + d_iY 
    \end{cases}
    \]
    with $b_i \neq d_i$ for all $i$. Let $(x_1,y_1), \cdots , (x_m,y_m)$ be generic solutions of any of these systems. 
    
    If $\mathrm{trdeg}(x_1,y_1, \cdots , x_m,y_m/\mathbb{C}) < 2m$, then either:
    \begin{itemize}
        \item $(b_i,d_i) = (b_j,d_j)$ for some $i ,j$ and $(c_i x_i, a_i y_i) = (c_j x_j, a_j y_j)$,
        \item $(b_i,d_i) = (d_j,b_j)$ for some $i , j$ and $(c_ix_i,a_iy_i) = (a_jy_j, c_jx_j)$.
    \end{itemize}
\end{thm}

\begin{proof}
    Note that we may as well assume, maybe by repeating and/or removing some of the systems, that $n=m$ and $(x_i,y_i)$ is a solution of $LV_{a_i,b_i,c_i,d_i}$ for all $i$. We do so in the rest of the proof. 
    
    Let $K$ be the subfield of $\mathbb{C}$ generated by the coefficients of all the $LV_{a_i,b_i,c_i,d_i}$. The fact that $\mathrm{trdeg}(x_1,y_1, \cdots , x_m,y_m/\mathbb{C})$ is strictly less than $2m$ implies that there is some countable field $F$ with $K<F<\mathbb{C}$ such that $x_1 , y_1 , \cdots x_m, y_m$ has transcendence degree strictly less than $2m$ over $F$

    Note that for all $i$, the type $\tp(x_i,y_i/\mathbb{Q}(a_i,b_i,c_i,d_i))$ is orthogonal to $\mathbb{C}$, and in particular $x_i y_i \forkindep_{\mathbb{Q}(a_i,b_i,c_i,d_i)} F$. Therefore, all the types $\tp(x_i,y_i/F)$ are minimal and disintegrated. Pick a maximal independent subset of the $(x_i,y_i)$, which we may assume, by reordering, to be of the form $(x_1,y_1), \cdots , (x_{j-1},y_{j-1})$ for some $j < n$. So we have both 
    \[(x_{j},y_{j}) \not\forkindep_F (x_1,y_1), \cdots ,(x_{j-1},y_{j-1})\]
    and
    \[(x_1,y_1) \forkindep_F (x_2,y_2), \cdots , (x_{j-1},y_{j-1}) \text{ .}\]
    By \cite[Chapter XVI, Lemma 2.4]{baldwin2017fundamentals} (keeping in mind that a disintegrated type is trivial in the sense defined by Baldwin in that book), we find that $(x_{j},y_j)$ forks over $F$ with either $(x_1,y_1)$ or $(x_2,y_2), \cdots , (x_{j-1},y_{j-1})$. Iterating this process, we obtain that $\mathrm{trdeg}(x_i,y_i,x_j,y_j/F) < 4$ for some $i \neq j$. In other words, we can conclude that $\mathrm{trdeg}(x_i,y_i,x_j,y_j/F) < 4$. Applying rescaling transformations, we obtain $\mathrm{trdeg}(c_ix_i,a_iy_i,c_jx_j,a_jy_j/F) < 4$. Since $(c_ix_i,a_iy_i)$ and $(c_jx_j,a_j y_j)$ are solutions of $LV_{1,b_i,1,d_i}$ and $LV_{1,b_j,1,d_j}$, respectively, Theorem \ref{theo: ortho-desin} implies that either $(b_i,d_i) = (b_j,d_j)$ and $(c_ix_i , a_i y_i) = (c_j x_j , a_j d_j)$, or $(b_i,d_i) = (d_j, b_j)$ and $(c_i x_i , a_i y_i) = (a_j y_j, c_j x_j)$.
\end{proof}

\subsection{The \texorpdfstring{$b=d$}{b=d} case}\label{subsec: b is d}

In this subsection, we deal with the special, non-strongly minimal case of $LV_{a,b,c,d}$ where $b = d$ by computing all possible irreducible invariant algebraic curves. Note that for this specific case, we have a simple, explicit expression for the solutions, due to Varma \cite{varma1977exact}. 
The set of solutions is given by:
\[
\begin{cases}
    x = \frac{\alpha e^{bt}}{a-ce^{\frac{\alpha e^{bt}- \beta}{b}}} \\
    y = \frac{\alpha e^{bt}}{c e^{-\frac{\alpha e^{bt}- \beta}{b}} - a}
\end{cases}
\]
\noindent where $\alpha, \beta$ are arbitrary complex numbers.

As previously discussed, we can reduce to studying the system $LV_{1,b,1,b}$
\[
\begin{cases}
    X' = XY + bX \\
    Y' = XY + bY
\end{cases}
\]
for some complex number $b$. Let $s$ be the associated vector field on the affine plane.

By our previous work, we know that there are three complex invariant, irreducible algebraic curves, given by $X=0$, $Y=0$ and $X-Y=0$. In fact, one can compute that $\mathcal{L}_s(X-Y) = b(X-Y)$. This implies the following, which was already observed in \cite[subsubsection 3.1.1]{duan2025algebraic}:

\begin{prop}
    Let $z$ be any solution of $Z' = bZ$. The polynomial $X-Y-z$ is invariant over $\mathbb{C}(z)$.
\end{prop}

\begin{proof}
    We compute:
    \[\mathcal{L}_s(X-Y-z) = XY+bX-XY-bY-bz = b(X-Y-z) \text{ .}\]
\end{proof}

We show that these are the only invariant irreducible curves. To do so, we use the model-theoretic notions of weight and domination. We give a brief explanation in the following paragraphs.

Let $p$ be a type over any field $F$, and assume that for some (any) realization $a$ of $p$, the differential field $F\langle a \rangle$ has finite transcendence degree over $F$. The \emph{preweight} of $p$, denoted $\mathrm{prwt}(p)$, is the largest $N \in \mathbb{N}$ such that there exist a realization $a$ of $p$ and an $F$-independent set $b_1, \cdots, b_N$ such that $a$ forks with $b_i$ over $F$ for all $1 \leq i \leq N$. The \emph{weight} of $p$, denoted $\mathrm{wt}(p)$, is the supremum, over all non-forking extensions $q$ of $p$, of $\mathrm{prwt}(q)$. Note that both weight and preweight are bounded by the order of $p$, and thus are well-defined natural numbers for finite order types (they are also well-defined for infinite order types, but this requires more model theoretic machinery). 

In what follows, we will be almost exclusively be concerned with types of weight $1$. If $\mathrm{wt}(p)=1$, then the weight of $p$ is equal to the preweight. Concretely, if $\mathrm{wt}(p)=1$, then for all $a$ realising $p$ and $b_1,b_2$ such that $a$ forks with $b_i$ over $F$ for $i=1,2$, we have that $b_1$ and $b_2$ fork over $F$. 

One way that two types can be related is for $p$ to `control' the possible differential-algebraic relations of $q$. Let $p\in S(K)$ and $q\in S(F)$ be two types with $F$ and $K$ algebraically closed differential fields. We say $p$ \emph{dominates} $q$ if there is some differential field $L$ containing $F$ and $K$, some $a\models p{\upharpoonright}_L$ and $b\models q{\upharpoonright}_L$ such that for every set $D$, if $a\forkindep_LD$ then $b\forkindep_LD.$ We say that $p$ and $q$ are \emph{domination-equivalent} if $p$ dominates $q$ and $q$ dominates $p$. If $p$ and $q$ are both weight 1 types, then $p$ is domination-equivalent to $q$ if and only if $p$ is non-orthogonal to $q$. As an example, minimal types have weight 1. 

When $p$ is a finite rank type in $\mathrm{DCF}_0$, $p$ has finite weight and is domination equivalent to $r_1\otimes\cdots\otimes r_n$ where $r_1,...,r_n$ are minimal types, potentially over a differential extension of $K$, and $n=\mathrm{wt}(p)$. See \cite[Chapter 1, sections 4.2 and 4.4]{pillay1996geometric} for a discussion of the above facts.

We will now compute the weight of the type of interest. 

\begin{lem}\label{lem: weight}
    Let $p$ be the generic type of $LV_{1,b,1,b}$ over $F = \mathbb{Q}(b)^{\mathrm{alg}}$. Then $\mathrm{wt}(p) = 1$.
\end{lem}

\begin{proof}
    By \cite[Theorem 4.7]{eagles2024internality}, the type $p$ is two-step analysable in the constants, but not almost internal to the constants. It follows that $p$ is domination-equivalent to the generic type of the constants, see \cite[Theorem 3.11]{eagles2025domination}. This implies that $\mathrm{wt}(p) = 1$.
\end{proof}

We can now prove: 

\begin{thm}\label{theo: LV b=d}
    The only irreducible invariant curves of $LV_{a,b,c,b}$ are given by $X=0$, $Y=0$, $cX-aY=0$ (the complex ones) and $cX-aY -z$ for $z$ a solution of $Z' = bZ$.
\end{thm}

\begin{proof}
    We assume, for now, that $a = c = 1$. Let $F = \mathbb{Q}(b)^{\mathrm{alg}}$.

    We have already classified the complex invariant irreducible curves in Lemma \ref{lem: complex-inv-curves}, so we only need to take care of those not defined over $\mathbb{C}$.

    Let $C$ be an invariant irreducible curve, defined over some differential field $F < K$ not contained in $\mathbb{C}$, and consider some generic, non-degenerate, point $(u,v)$ in $C$, which exists as $C$ is not defined over $\mathbb{C}$. We may take $K = F\langle w \rangle$ for some tuple $w$. The transcendence degree of $(u,v)$ over $K$ is $1$, so $(u,v)$ forks with $w$ over $F$. Note that we also have that $(u,v)$ forks with $u-v$ over $F$. Since $\mathrm{wt}(p) = 1$, this implies that $u-v$ forks with $w$ over $F$. Because $(u-v)' = b(u-v)$, we have $\mathrm{trdeg}(u-v/F) = 1$, so we conclude that $u-v \in K^{\mathrm{alg}}$.

    Now, the type of $(u,v)$ over $F\langle u-v \rangle$ contains the formula $x-y-(u-v) = 0$, i.e. $(u,v)$ is in the invariant curve given by $X-Y - (u-v)$. Therefore, the type of $(u,v)$ over $K$ implies that $(u,v)$ is contained in finitely many irreducible invariant algebraic curves of the form $X-Y-z$, for some $z$ solution of $Z'=bZ$. By irreducibility, the curve $C$ has to be equal to one of them. Keeping in mind the transformation applied at the start, this yields the theorem.
\end{proof}

\bibliography{biblio}
\bibliographystyle{plain}

\end{document}